\documentclass[12pt]{article}

 \usepackage{amsmath}
\usepackage{amsfonts}
\usepackage{color}
\usepackage{bm}

\def\prob {{\sf Pr}}
\def\sg {{\sf g}}

\def\sm {{\sf p}}

\newcommand{\cA}{{\mathfrak A}}

\newcommand{\cQ}{{\mathfrak Q}}

\newcommand{\cM}{{\mathfrak M}}

\newcommand{\cR}{{\mathfrak R}}

\newcommand{\cv}{{\mathfrak v}}

\newtheorem{proposition}{Proposition}[section]
\newtheorem{example}{Example}[section]
\newtheorem{remark}{Remark}[section]
\newtheorem{definition}{Definition}[section]

\newtheorem{assumption}{Assumption}[section]

\newcommand{\R}{{\cal R}}

\def\argmin{\mathop{\rm arg\,min}}
\newcommand{\var}{{\rm Var}}

\newcommand{\Q}{{\cal Q}}
\newcommand{\U}{{\cal U}}

\newcommand{\Z}{{\cal Z}}
\newcommand{\F}{{\cal F}}

\newcommand{\G}{{\cal G}}
\newcommand{\I}{{\cal I}}
\newcommand{\J}{{\cal J}}

\newcommand{\X}{{\cal X}}

\newcommand{\be}{\begin{equation}}
\newcommand{\ee}{\end{equation}}

\def\dist{\mathop{\rm dist}}
\def\w{\omega}
\def\e{\epsilon}
\def\O{\Omega}
\def\vv{\vartheta}
\def\e{\varepsilon}

\def\dom{{\rm dom}}

\def\Gr {{\sf Gr}}

\def\und {\underline{\cal Q}}
\def\undv {\underline{V}}

\def\argmin{\mathop{\rm arg\,min}}
\def\argmax{\mathop{\rm arg\,max}}

\newcommand{\avr}{{\sf AV@R}}

\def\bbr{{\Bbb{R}}} 
\def\bbe{{\Bbb{E}}} 

\def\bbp{{\Bbb{P}}}

\def\eqnok#1{(\ref{#1})}
\newcommand{\tsum}{{\textstyle\sum}}
\def\argmin{\mathop{\rm arg\,min}}
\def\argmax{\mathop{\rm arg\,max}}

\newcommand{\beq}{\begin{equation}}
\newcommand{\eeq}{\end{equation}}

\def\dist{\mathop{\rm dist}}
\def\gap{\mathop{\rm gap}}

\begin{document}

\begin{titlepage}
\title{\bf Numerical Methods for Convex Multistage Stochastic Optimization}
\author{ \bf   Guanghui Lan\thanks{Georgia Institute of Technology, Atlanta, Georgia
30332, USA,
\tt{george.lan@isye.gatech.edu}\newline Research of this author was
partially supported by the NSF grant
DMS-1953199 and the NSF AI Institute grant NSF-2112533.}\and { \bf Alexander Shapiro}
\thanks{Georgia Institute of Technology, Atlanta, Georgia
30332, USA, \tt{ashapiro@isye.gatech.edu}\newline
Research of this
author was partially supported by Air Force Office of Scientific Research (AFOSR)
under Grant FA9550-22-1-0244.}}


\maketitle
\begin{abstract}

Optimization problems involving sequential decisions in  a  stochastic environment    were studied  in  Stochastic Programming (SP), Stochastic Optimal Control  (SOC) and Markov Decision Processes (MDP). In this paper we mainly concentrate on SP and  SOC modelling   approaches. In these frameworks there are natural situations  when the considered problems are  convex.
Classical approach to sequential optimization is based on dynamic programming. It has the problem of the so-called ``Curse of Dimensionality", in that
its computational complexity increases exponentially with increase of dimension of
state variables. Recent progress in solving convex multistage   stochastic  problems is based on cutting planes
approximations of the cost-to-go (value) functions of dynamic programming equations. Cutting planes type
algorithms in dynamical settings is one of the main topics of this paper. We also discuss Stochastic Approximation type methods applied to multistage stochastic optimization problems.
From the computational complexity point of view, these two types of methods seem
to be complimentary to each other. Cutting plane type methods can
handle multistage problems with a large number of stages, but a relatively smaller number
of state (decision) variables. On the other hand, stochastic approximation type methods can
only deal with a small number of stages, but a large number of decision variables.
\end{abstract}
\end{titlepage}

\tableofcontents

\setcounter{equation}{0}
\section{Introduction}
\label{sec-intr}

 Traditionally different communities  of researchers dealt with optimization problems involving uncertainty, modelled in stochastic terms, using different terminology and modeling frameworks. In this respect we can point to the fields
 of Stochastic Programming (SP), Stochastic Optimal Control  (SOC) and Markov Decision Processes (MDP). Historically the developments in SP  on one hand  and SOC and MDP on the other,   went along different directions with different modeling frameworks  and solution methods. In this paper
 we mainly concentrate on SP and  SOC
   approaches. In these frameworks there are natural situations  when the considered problems are {\em convex}. In continuous optimization convexity is the main consideration  allowing development of efficient  numerical algorithms, \cite{Nesnem}.

    The main goal of this paper is to present some recent  developments in numerical approaches to
   solving {\em convex} optimization problems involving sequential decision making. We do not try   to give a comprehensive  review of the subject with a complete list of references, etc.
{\em Rather the aim is to present  a certain point of view about  some recent developments in solving    {\em convex }  multistage stochastic problems.}

  The SOC modeling is convenient for demonstration of some basic ideas since   on one hand it can be naturally written in the MDP terms, and on the other hand it can be formulated in the Stochastic Programming framework.
Stochastic Programming (SP) has a long history.  Two stage stochastic programming  (with recourse) was introduced in Dantzig \cite{dant:55} and
Beale \cite{beal:55},  and was intrinsically connected with linear programming. From the beginning SP was aimed at numerical solutions.  Until about twenty  years ago, the modeling approach to two and multistage SP was predominately based   on construction of scenarios represented by scenario trees. This approach  allows to formulate the so-called deterministic equivalent optimization problem with the number of decision   variables more or less proportional to the number of scenarios.
When the deterministic equivalent could be represented as  a linear program, such problems
were considered to be numerically solvable. Because of that,   the topic of SP   was often  viewed as a large scale linear programming. We can refer for presentation of
these developments to Birge \cite{birge2011} and references there in.

From the point of view of the scenarios  construction approach  there is no much difference between   two stage and multistage SP.  In both cases the numerical effort in solving the deterministic equivalent  is more or less proportional to the number of generated scenarios. This view on SP  started to change with developments of randomization methods and the sample complexity theory.
From the     point of view  of solving  deterministic equivalent, even two stage linear stochastic programs are computationally intractable; their computational complexity is   \#P-hard for a sufficiently high
accuracy
(cf., \cite{dyer},\cite{HKW}). On the other hand,
under reasonable assumptions,  the number   of randomly generated scenarios  (by Monte Carlo sampling techniques),  which are required to solve two stage SP problems  with accuracy $\e>0$ and high probability is of order $O(\e^{-2})$,  \cite[Section 5.3]{SDR}. While randomization methods were reasonably successful in solving two stage problems, as far as multistage SP is concerned the situation is different. Number of scenarios needed to solve multistage SP problems grows exponentially with increase of the number of stages, \cite{shanem},\cite[Section 5.8.2]{SDR}.

Classical approach to sequential optimization is based on dynamic programming, \cite{Bel}. Dynamic programming  also has a long history and is in the heart of the SOC and MDP modeling. It has the problem of the so-called  ``Curse of Dimensionality", the term coined by Bellman.
Its  computational complexity increases exponentially with increase of dimension of state variables. There is a large literature trying to deal with this problem by using various approximations of dynamic programming equations (e.g., \cite{powell}).  Most of these methods are heuristics and do not give verifiable guarantees  for  the accuracy of obtained solutions.

Recent progress in solving convex multistage SP problems is based on cutting planes approximations of the cost-to-go (value) functions of dynamic programming equations.
These methods allow to give an estimate of the error of  the computed solution.
Cutting planes type algorithms  in dynamical settings is one of  the main topics of this paper. We discuss such algorithms in the frameworks of SP and SOC.
Moreover, we will also present extensions of stochastic approximation (a.k.a. stochastic gradient descent) type methods~\cite{RobMon51-1,nemyud:83,NJLS09-1,LanBook2020}
for multistage stochastic optimization, referred to as dynamic stochastic approximation in \cite{LanZhou17-1}.
From the computational complexity point of view, these two types of methods seem to be complimentary to each other in the following sense.
Cutting plane type methods can handle multistage problems with a large number of stages, but a relatively
smaller number of state (decision) variables. On the other hand, stochastic approximation type methods can only deal with a small number of stages, but
a large number of decision variables. These methods share the following common features:
(a) both methods utilize the convex structure of the cost-to-go (value) functions of dynamic programming equations,
(b) both methods do not require explicit discretization of state space,
(c) both methods guarantee the convergence to the global optimality,
 (d) rates of convergence for both methods have been established.
\\

 We use the following notation and terminology  throughout the paper. For $a\in \bbr$ we denote $[a]_+:=\max\{0,a\}$. Unless stated otherwise $\|\cdot\|$ denotes Euclidean norm in $\bbr^n$.
  By $\dist (x, S):=\inf_{y\in S}\|x-y\|$  we denote the distance from  a point $x\in\bbr^n$ to a set $S\subset \bbr^n$. 
   We   write $x^\top y$ or $\langle x,y\rangle$ for the scalar product 
   $\sum_{i=1}^n x_i y_i$ 
 of  vectors $x,y\in \bbr^n$.
  It is said that a set $S\subset \bbr^n$ is polyhedral if it can be represented by a finite number of affine constraints, it is said that a
  function $f:\bbr^n\to \bbr$ is polyhedral if it can be represented as maximum of a finite number of affine functions.   For
a process $\xi_1,\xi_2,...$,  we denote by $\xi_{[t]}=(\xi_1,...,\xi_t)$ its history up to time $t$. By $\bbe_{|X}[\,\cdot\,]$ we denote the conditional expectation,  conditional on  random variable (random vector)   $X$.  We use the same notation  $\xi_t$ viewed as a random vector or as a vector variable, the  particular meaning will be clear from the context. 
For a probability space $(\O,\F,\bbp)$, 
by $L_p(\O,\F,\bbp)$, $p\in [1,\infty)$,  we denote the space of random variables $Z:\O\to \bbr$ having finite $p$-th order moment, i.e., such that $\int |Z|^p d\bbp <\infty$. Equipped with norm $\|Z\|_p:=
\left(\int  |Z|^p d\bbp\right)^{1/p}$,  $L_p(\O,\F,\bbp)$ becomes a Banach space.
The dual of $\Z:=L_p(\O,\F,\bbp)$ is the space $\Z^*=L_q(\O,\F,\bbp)$ with $q\in (1,\infty]$   such that $1/p+1/q=1$.  

 \setcounter{equation}{0}
\section{Stochastic Programming}
\label{sec-sp}

 In Stochastic Programming (SP) the $T$-stage  optimization problem  can be written as
 \begin{eqnarray}
 \label{stpr-1}
   & \min\limits_{x_t\in \X_t}& \bbe\Big[\sum_{t=1}^T f_t(x_t,\xi_t)\Big]\\
   \label{stpr-1a}
 & {\rm s.t}&  B_t x_{t-1} +A_t x_t=b_t,\;t=1,...,T.
 \end{eqnarray}
Here $f_t:\bbr^{n_t}\times \bbr^{d_t}\to \bbr$ are objective functions,
$\xi_t\in \bbr^{d_t}$  are random vectors, $\X_t\subset \bbr^{n_t}$ are nonempty closed sets,  $B_t=B_t(\xi_t)$,  $A_t=A_t(\xi_t)$ and $b_t=b_t(\xi_t)$ are   matrix and vector  functions of $\xi_t$,  $t=1,...,T $, with $B_1=0$.

Constraints  \eqref{stpr-1a} often represent some type of   balance equations  between successive stages, while the set  $\X_t$ can be viewed as representing  local constraints at time $t$. It is possible to consider balance equations in a more general form, nevertheless formulation \eqref{stpr-1a} will be sufficient for our purpose of discussing convex problems. In particular, if the objective functions  $f_t(x_t,\xi_t):=c_t^\top x_t$ are linear, with $c_t=c_t(\xi_t)$, and the sets $\X_t$  are polyhedral,   then problem \eqref{stpr-1} - \eqref{stpr-1a} becomes a linear multistage stochastic program.

The sequence of random  vectors $\xi_1,...$,  is   viewed as a data process.
Unless stated otherwise we make the following assumption throughout the paper.
\begin{assumption}
  \label{assum-ind}
  Probability distribution of the random process $\xi_t$, $t=1,...$,   does not depend on our decisions, i.e., is independent of the chosen decision policy.
\end{assumption}

At every time period $t=1,...,$  we have information about (observed)   realization   $\xi_{[t]}=(\xi_1,...,\xi_t)$  of the data process. Naturally our decisions  should be based on the information available at time of the decision and should not depend on the future unknown values $\xi_{t+1},...,\xi_{T}$;
this is the so-called principle of {\em nonantisipativity}.
That is, the decision variables $x_t=x_t(\xi_{[t]})$ are {\em functions} of the data process, and a sequence of such functions for every stage $t=1,...,T$,  is called a {\em policy} or a decision rule.
The optimization in \eqref{stpr-1} - \eqref{stpr-1a} is performed over  policies satisfying the feasibility constraints.
  The feasibility constraints of problem \eqref{stpr-1} - \eqref{stpr-1a} should be satisfied with probability one, and the expectation is taken with respect to the probability distribution of the random vector $\xi_{[T]}=(\xi_1,...,\xi_T)$.
It could be noted that the dependence of decision variables  $x_t$,  as well as parameters $(B_t$,$A_t$,$b_t)$, on the data process  is often suppressed in the notation.
It is worthwhile to emphasize that it suffices to consider policies depending only   on the history of the data process without  history of the decisions, because of Assumption \ref{assum-ind}.

The first stage decision $x_1$ is of the main interest and assumed to be
 {\em deterministic}, i.e., made before observing future realization of the data process. In that framework vector $\xi_1$ and the corresponding parameters   $(A_1,b_1)$  are
 deterministic, i.e., their  values are  known at the time of  first stage decision. Also  the first stage objective function  $f_1(x_1)=f_1(x_1,\xi_1)$ is a function of $x_1$ alone.

We can write the following dynamic programming equations   for problem \eqref{stpr-1} - \eqref{stpr-1a}  (e.g., \cite[Section 3.1]{SDR}). At the last stage the  cost-to-go (value) function is
\begin{equation}\label{dynpr-1}
Q_T(x_{T-1},\xi_T)=  \inf_{x_T\in \X_T} \left\{ f_T(x_T,\xi_T): B_T x_{T-1} +A_T x_T=b_T\right\}
\end{equation}
and then going backward in time for $t=T-1,...,2$,  the cost-to-go  function is
\begin{equation}\label{dynpr-2}
Q_t(x_{t-1},\xi_{[t]})=  \inf_{x_t\in \X_t} \left\{ f_t(x_t,\xi_t)+\bbe_{|\xi_{[t]}}[Q_{t+1}(x_t,\xi_{[t+1]})]: B_t  x_{t-1} +A_t x_t=b_t\right\}.
\end{equation}
The optimal value  of problem \eqref{stpr-1} - \eqref{stpr-1a} is  given by the optimal value of the first stage problem
\begin{equation}\label{dynpr-3}
  \min\limits_{x_1\in \X_1}  f_1(x_1)+ \bbe\left[  Q_2(x_1,\xi_2)\right]\;
  {\rm s.t.}\;  A_1x_1=b_1.
 \end{equation}

The optimization in \eqref{dynpr-1} - \eqref{dynpr-2} is over (nonrandom)  variables $x_t\in \bbr^{n_t}$    (recall that $B_t=B_t(\xi_t)$,  $A_t=A_t(\xi_t)$ and $b_t=b_t(\xi_t)$  are    functions of $\xi_t$). The corresponding optimal policy is defined by $\bar{x}_t=\bar{x}_t(\xi_{[t]})$, $t=2,...,T$, where
\begin{equation}\label{dynpr-4}
 \bar{x}_t \in   \argmin_{x_t\in \X_t} \left\{ f_t(x_t,\xi_t)+\bbe_{|\xi_{[t]}}[Q_{t+1}(x_t,\xi_{[t+1]})]: B_t \bar{x}_{t-1} +A_t x_t=b_t\right\},
\end{equation}
with the expectation term of  $Q_{T+1}(\cdot,\cdot)$  at the last stage omitted.  Of the main interest is the first stage solution given by the optimal solution $\bar{x}_1$  of the first stage problem   \eqref{dynpr-3}.

Since    $\bar{x}_{t-1}$ is a function of $\xi_{[t-1]}$,
  the cost-to-go  function $Q_t(\bar{x}_{t-1},\xi_{[t]})$
 can be considered as a function of  $\xi_{[t]}$. It
represents the optimal value of problem
\begin{equation}
 \label{sptail}
    \min\limits_{x_\tau\in \X_\tau}  \bbe\Big[\sum_{\tau=t+1}^T f_\tau(x_\tau,\xi_\tau)\Big]\; \;
     {\rm s.t.}\;  B_\tau x_{\tau-1} +A_\tau x_\tau=b_\tau,\;\tau=t+1,...,T,
 \end{equation}
conditional on realization $\xi_{[t]}$ of the data process.
 This is the dynamic programming (Bellman) principle of optimality. This is also the motivation for the name ``cost-to-go" function.  Another name for  the cost-to-go functions is ``value functions".

\begin{remark}
\label{rem-conv-sp}
{\rm
Let us recall the following simple result. Consider an extended real valued  function $\phi:\bbr^n\times \bbr^m\to \bbr\cup\{+\infty\}$. We have that if $\phi(x,y)$ is a convex function of $(x,y)\in \bbr^n\times \bbr^m$, then the min-value function  $\varphi(x):=\inf_{y\in \bbr^m}\phi(x,y)$ is also convex. This can be applied to establish convexity of the cost-to-go functions. That is, suppose that   the function $f_T(x_T,\xi_T)$ is  convex in $x_T$ and the set  $\X_T$ is  convex. Then the extended real valued function taking value $f_T(x_T,\xi_T)$ if $x_T\in \X_T$ and $B_T x_{T-1} +A_T x_T=b_T$,  and value $+\infty$ otherwise, is convex. Since   the cost-to-go function
$Q_T(x_{T-1},\xi_T)$    can be viewed as the min-function of
this extended real valued function, it follows that it is convex in $x_{T-1}$. Also note that if the cost-to-go function $Q_{t+1}(x_t,\xi_{[t+1]})$  is convex in $x_t$, then its  expected value in  \eqref{dynpr-2} is also convex. Consequently by induction going backward in time it is not difficult to show the following convexity property of the  cost-to-go functions.
\begin{proposition}
 \label{pr-spconvex}
 If the objective functions $f_t(x_t,\xi_t)$ are convex in $x_t$ and the sets $\X_t$ are convex, $t=1,...,T$,  then the cost-to-go functions $Q_{t}(x_{t-1},\xi_{[t]})$, $t=2,...,T$,  are convex in $x_{t-1}$ and the first stage problem \eqref{dynpr-3} is convex. In particular such convexity   follows   for  {\em linear} multistage stochastic programs.
\end{proposition}

Such convexity property is crucial for development   of efficient numerical algorithms. It is worthwhile to note that Assumption \ref{assum-ind} is essential here. Dependence of the distribution of the random data on decisions typically destroys convexity even in the linear case.
} $\hfill \square$
\end{remark}

The dynamic programming equations reduce the problem from optimization over policies, i.e., over  infinite dimensional functional spaces,  to evaluation of  the cost-to-go functions $Q_t(x_{t-1},\xi_{[t]})$    of finite dimensional vectors. However, dependence of the cost-to-go functions on the whole  history of the data process makes it numerically unmanageable with increase of the number of stages. The situation is simplified dramatically if  we make the following assumption of {\em stagewise independence}.

\begin{assumption} [Stagewise independence]
  \label{assum-stage}
  Probability distribution of   $\xi_{t+1}$   is independent of
  $\xi_{[t]}$ for $t\ge 1$.
\end{assumption}

Under   Assumption \ref{assum-stage} of stagewise independence, the conditional expectations in dynamic equations   \eqref{dynpr-2}  become the  corresponding unconditional expectations. It is straightforward to show then by induction,  going backward in time, that the dynamic equations  \eqref{dynpr-2} can be written as
\begin{equation}\label{dynind-1}
Q_t(x_{t-1},\xi_{t})=  \inf_{x_t\in \X_t} \left\{ f_t(x_t,\xi_t)+\Q_{t+1}(x_{t}): B_t x_{t-1} +A_t x_t=b_t\right\},
\end{equation}
where
\begin{equation}\label{dynind-1a}
  \Q_{t+1}(x_{t}):=  \bbe [Q_{t+1}(x_t,\xi_{t+1})],
\end{equation}
with the expectation taken with respect to the marginal distribution of $\xi_{t+1}$.
In that setting one needs  to keep track of the (expectation) cost-to-go functions $\Q_{t+1}(x_{t})$ only. The optimal policy is defined by the equations
\begin{equation}\label{dynind-2}
 \bar{x}_t \in   \argmin_{x_t\in \X_t}   \left\{ f_t(x_t,\xi_t)+\Q_{t+1}(x_{t}): B_t \bar{x}_{t-1} +A_t x_t=b_t\right\}.
\end{equation}
Note that here the optimal policy values can be computed iteratively, starting from optimal solution $\bar{x}_1$ of the first stage problem, and then   by sequentially  computing a minimizer in the right hand side of \eqref{dynind-2} going forward in time for  $t=2,...,T$. Of course,  this procedure requires availability of the cost-to-go functions $\Q_{t+1}(x_{t})$.

It is worthwhile to note that here the  optimal policy decision $\bar{x}_t$  can be viewed as a function of $\bar{x}_{t-1}$ and $\xi_t$, for $t=2,...,T$. Of course,
$\bar{x}_{t-1}$ is a function of $\bar{x}_{t-2}$ and $\xi_{t-1}$, and so on. So eventually $\bar{x}_{t}$ can be represented as a function of the history $\xi_{[t]}$ of the data process.    Assumption \ref{assum-ind} was crucial for  this conclusion, since then we   need to consider only the history of the data process rather than also the   history of our decisions. This is an important point, we will discuss this later.

\begin{remark}
\label{rem-mark}
{\rm
In SP setting the process $\xi_t$ is viewed as a data process, and the assumption of stagewise independence could be unrealistic.  In order to model stagewise dependence we assume  the following Markovian property of the data process: the conditional distribution of $\xi_{t+1}$ given $\xi_{[t]}$ is the same as  the conditional distribution of $\xi_{t+1}$ given $\xi_{t}$. In such Markovian setting the dynamic equations \eqref{dynpr-2}  become
\begin{equation}\label{dynind-3}
Q_t(x_{t-1},\xi_{t})=  \inf_{x_t\in \X_t} \left\{ f_t(x_t,\xi_t)+\bbe_{|\xi_{t}}[Q_{t+1}(x_t,\xi_{t+1})]: B_t  x_{t-1} +A_t x_t=b_t\right\}.
\end{equation}
 That is,   compared with the stagewise independent case, in the Markovian setting   the expected cost-to-go function
 \begin{equation}\label{dynind-3a}
 \Q_{t+1}(x_{t},\xi_t)=  \bbe_{|\xi_{t}}[Q_{t+1}(x_t,\xi_{t+1})]
 \end{equation}
also depends on $\xi_t$.
} $\hfill \square$
\end{remark}

\setcounter{equation}{0}
\section{Stochastic Optimal Control}
\label{sec-soc}

Consider the classical Stochastic Optimal Control  (SOC)  (discrete time, finite horizon) model (e.g. , Bertsekas and   Shreve \cite{ber78}):
\begin{eqnarray}
\label{soc}
&\min\limits_{u_t\in \U_t(x_t)} &\bbe^\pi\Big [ \sum_{t=1}^{T}
c_t(x_t,u_t,\xi_t)+c_{T+1}(x_{T+1})
\Big],\\
&{\rm s.t.} &
\label{soc-b}
x_{t+1}=F_t(x_t,u_t,\xi_t),\;
t=1,...,T.
\end{eqnarray}
Variables  $x_t\in \bbr^{n_t}$, $t=1,...,T+1$, represent the state  of the system,   $u_t\in \bbr^{m_t}$,  $t=1,...,T$, are controls,   $\xi_t\in \bbr^{d_t}$, $t=1,...,T$, are random vectors, $c_t:\bbr^{n_t}\times\bbr^{m_t}\times\bbr^{d_t}\to \bbr$, $t=1,...,T$, are cost functions, $c_{T+1}(x_{T+1})$ is a final cost function,  $F_t:\bbr^{n_t}\times\bbr^{m_t}\times\bbr^{d_t}\to \bbr^{n_{t+1}}$ are (measurable) mappings and $\U_t(x_t)$  is a (measurable)  multifunction mapping  $x_t\in \bbr^{n_t}$ to a subset of   $\bbr^{m_t}$.
Values  $x_1$  and $\xi_0$ are  deterministic  (initial conditions); it is also possible to view $x_1$ as random with a given distribution, this is not essential for the following discussion.
The optimization in \eqref{soc} - \eqref{soc-b} is performed over policies $\pi$  determined by  decisions  $u_t$ and state variables $x_t$.
This is emphasized in the notation $\bbe^\pi$, we are going to discuss this below.


Unless stated otherwise we assume   that probability distribution of the random process $\xi_t$  does not depend on our decisions
(Assumption \ref{assum-ind}).
We can consider  problem \eqref{soc}-\eqref{soc-b} in the framework of SP   if we view $y_t=(x_{t},u_t)$, $t=1,...,T$,  as decision variables.
 Suppose  that the   mapping  $F_t(x_t,u_t,\xi_t)$ is affine, i.e.,
\begin{equation}\label{affine}
 F_t(x_t,u_t,\xi_t):=A_t x_t +B_t u_t +b_t, \;t=1,...,T,
\end{equation}
where $A_t=A_t(\xi_t) $, $B_t=B_t(\xi_t)$ and $b_t =b_t(\xi_t)$ are matrix and vector valued functions of $\xi_t$.
The constraints $u_t\in \U_t(x_t)$ can be viewed as local constraints in the SP framework; in particular if
$\U_t(x_t)\equiv \U_t$ are independent of $x_t$, then $\U_t$ can be viewed as a counterpart  of the set $\X_t$ in problem \eqref{stpr-1} - \eqref{stpr-1a}.
For the affine mapping of the form \eqref{affine}, the state equations \eqref{soc-b} become
 \begin{equation}\label{state-eq}
 x_{t+1}-A_t x_t -B_t u_t=b_t,\;
t=1,...,T.
 \end{equation}
These state equations are linear in $y_t=(x_{t},u_t)$ and can be viewed as a particular case of  the balance equations \eqref{stpr-1a} of the stochastic program \eqref{stpr-1} - \eqref{stpr-1a}.
Note, however,     that here decision $u_t$ should be made before realization of $\xi_t$ becomes known, and  $x_{t}$ and $u_{t}$   are functions of $\xi_{[t-1]}$ rather than $\xi_{[t]}$.  We emphasize that  $x_{t}$ and $u_{t}$   can be considered as functions of the history $\xi_{[t-1]}$  of the  random process $\xi_t$,  without dependence on the history of the decisions.
As in the SP framework, this follows from the dynamic programming equations discussed below, and    Assumption \ref{assum-ind} is essential for that conclusion.

 There  is a  time shift in   equations \eqref{state-eq} as compared with the SP equations  \eqref{stpr-1a}.
 We emphasize that in both the SOC and SP approaches, the decisions are made based on  information available at time of the decision;
 this is the principle of nonanticipativity.
 Because of the time shift, the dynamic programming equations for the SOC problem are slightly different from the respective equations in  SP.  More important is that in the SOC framework there is a clear separation between the state and control variables.  This will have an important consequence  for  the respective dynamic programming equations which we consider next.

Similar to SP,
the dynamic programming equations for the SOC problem \eqref{soc} - \eqref{soc-b}
can be written as follows.
At the last stage, the value function
 $V_{T+1}(x_{T+1})=c_{T+1}(x_{T+1})$ and, going backward in time for $t=T,...,1$,  the value functions
 \begin{equation}\label{socvalfun-1}
   V_t(x_t,\xi_{[t-1]})=\inf\limits_{{u_t\in \U_t(x_t)}\atop{x_{t+1}=F_t(x_t,u_t,\xi_t)}
   }\bbe_{| \xi_{[t-1]}} \left [c_t(x_t,u_t,\xi_t)+
V_{t+1}\big(x_{t+1},\xi_{[t]} \big)\right].
 \end{equation}
 The optimization in the right hand side of \eqref{socvalfun-1}   is performed jointly in  $u_t$ and $x_{t+1}$.
 The optimal value of the SOC problem \eqref{soc}-\eqref{soc-b} is given by the first stage value function $V_1(x_1)$, and can be viewed  as a function of the initial conditions $x_1$.
 Note that because of the state equation \eqref{soc-b}, equations \eqref{socvalfun-1}
   can be written in the following equivalent form
 \begin{equation}
\label{soc-2}
V_t(x_t,\xi_{[t-1]})=\inf\limits_{u_t\in \U_t(x_t)}
\bbe_{| \xi_{[t-1]}} \left [c_t(x_t,u_t,\xi_t)+
V_{t+1}\big(F_t(x_t,u_t,\xi_t),\xi_{[t]} \big)\right].
\end{equation}

 The corresponding  optimal policy is defined by
  \begin{equation}\label{socvalfun-2}
  (\bar{u}_t,\bar{x}_{t+1})\in \argmin\limits_{{u_t\in \U_t(\bar{x}_t)}\atop{x_{t+1}=F_t(\bar{x}_t,u_t,\xi_t)}
   }\bbe_{| \xi_{[t-1]}} \left [c_t(\bar{x}_t,u_t,\xi_t)+
V_{t+1}\big(x_{t+1},\xi_{[t]} \big)\right].
 \end{equation}
  Similar to SP,    it follows that $\bar{u}_t$ can be considered as a function of $\xi_{[t-1]}$, and $\bar{x}_{t+1}$   as a function of $\xi_{[t]}$. That is,  $(\bar{u}_t,\bar{x}_{t})$  is a function of $\xi_{[t-1]}$, this was already mentioned in the previous paragraphs.
  Again Assumption \ref{assum-ind} is essential for the above derivations and conclusions.

 \begin{remark}\label{rem-convsoc}
 {\rm
 Let us discuss     convexity of the value functions. Consider graph of the  multifunction $\U_t(\cdot)$:
 \begin{equation}\label{graph}
   \Gr(\U_t):=\left \{ (x_t,u_t):u_t\in \U_t(x_t),\;x_t\in  \bbr^{n_t} \right\}.
 \end{equation}
 Note that minimization $\min_{u_t\in \U_t(x_t)}\varphi(x_t,u_t)$  of a function $\varphi:\bbr^{n_t}\times  \bbr^{m_t}\to \bbr\cup\{+\infty\}$ can equivalently be written as
  $\min_{u_t\in \U_t(x_t)}\hat{\varphi}(x_t,u_t)$,  where function $\hat{\varphi}(x_t,u_t)$ coincides with  $\varphi(x_t,u_t)$ when $u_t\in \U_t(x_t)$, and is $+\infty$ otherwise. Of course, the constraint $u_t\in \U_t(x_t)$ is the same as $(x_t,u_t)\in \Gr(\U_t)$. If function $\varphi(x_t,u_t)$ is convex and the set
$\Gr(\U_t)$ is convex, then the corresponding  function  $\hat{\varphi}(x_t,u_t)$ is convex. Consequently in that case  the min-function $\min_{u_t\in \U_t(x_t)}\varphi(x_t,u_t)$ is a convex function of $x_t$. By applying this observation and using induction   going backward in time, it is not difficult to show
that value functions
$V_{t}\big(x_{t},\xi_{[t-1]} \big)$ are convex in $x_{t}$, $t=1,...,T$, if the following assumption holds.

 \begin{assumption}
   \label{ass-conv}
   The function  $c_{T+1}(\cdot)$ is convex, and  for every $t=1,...,T,$  the cost function  $c_t(x_t,u_t,\xi_t)$  is convex in $(x_t,u_t)$, the mapping $F_t$ is affine of the form \eqref{affine}, the  graph $\Gr(\U_t)$ is  a convex subset of $\bbr^{n_t}\times  \bbr^{m_t}$, $t=1,...,T$.
 \end{assumption}
Note that, in particular, $\Gr(\U_t)$ is convex if $\U_t(x_t)\equiv \U_t$ is independent of $x_t$ and the set $\U_t$ is convex.
When  $\U_t(x_t)$ is dependent on $x_t$, we need some constructive way to represent it  in order to ensure convexity of $\Gr(\U_t)$. For example suppose that the feasibility sets are defined by inequality constraints as:
\begin{equation}\label{setrepr}
 \U_t(x_t):=\{u_t\in \bbr^{m_t}: g_{t}(x_t,u_t)\le 0\},
\end{equation}
 where $g_{t}:\bbr^{n_t}\times \bbr^{m_t} \to \bbr^{\ell_t}$. Then
 $\Gr(\U_t)=\{(x_t,u_t):g_t(x_t,u_t)\le 0\}$. Consequently if every component of
 mapping $g_t(x_t,u_t)$ is a convex function, then $\Gr(\U_t)$ is convex.
} $\hfill \square$
\end{remark}

In SOC framework  the   random process $\xi_1,...,$ is often viewed as a noise or disturbances (note that here $\xi_1$ is also random).
 In such cases the assumption of stagewise independence  (Assumption \ref{assum-stage}) is natural. In the stagewise independent case the dynamic programming equations \eqref{soc-2} simplify with the value functions $V_t(x_t)$, $t=1,...,T$,  being functions of the state variables only, that is
\begin{equation}
\label{soc-ind}
V_t(x_t)=\inf\limits_{u_t\in \U_t(x_t)}
\bbe \left [c_t(x_t,u_t,\xi_t)+
V_{t+1}\big(F_t(x_t,u_t,\xi_t) \big)\right],
\end{equation}
where the expectation is taken with respect to the marginal distribution of $\xi_t$.
Consequently the optimal policy is defined by the optimal controls   $\bar{u}_t=\pi_t(x_{t})$, where
\begin{equation}\label{conpol-2}
\bar{u}_t\in  \argmin_{u_t\in \U_t(x_t)}
\bbe  \left [c_t(x_t,u_t,\xi_t)+
V_{t+1}\big(F_t(x_t,u_t,\xi_t) \big)\right].
\end{equation}
 That is, in the stagewise independent case it is suffices to consider policies with controls of the form  $u_t=\pi_t(x_{t})$.

\begin{remark}
{\rm
Suppose   that  distribution of $\xi_t$  can  depend on    $(x_t,u_t)$, i.e., probability distribution $P_t(\cdot|x_t,u_t)$  of $\xi_t$ is conditional  on    $(x_t,u_t)$.  That is,  distribution of the random process is effected by our decisions and consequently  Assumption \ref{assum-ind}  is not satisfied.  Under    the stagewise independence Assumption \ref{assum-stage}, the dynamic programming equations \eqref{soc-ind} still hold and the optimal policy is determined by \eqref{conpol-2}.  However in that case $\bar{u}_t$ depends on
the history $(\bar{x}_1,\bar{u}_1,...,\bar{u}_{t-1}, \bar{x}_t)$ of the decision process. This history depends on the decisions (on the considered policy) and cannot be represented in terms of the process $\xi_t$ alone.  It also could be noted that in the stagewise independence case,  SOC can be formulated as Markov Decision Processes (MDP) with transition probabilities $P_t(\cdot|x_t,u_t)$   determined  by the state equations \eqref{soc-b}. We will not pursue the MDP formulation  in the following discussion, and will rather deal with the SOC model instead.}
$\hfill \square$
\end{remark}

As an example let us consider the classical inventory problem.

\begin{example}[Inventory problem]
\label{ex-inven}
{\rm
Consider problem  \cite{zipkin},
\begin{eqnarray}
\label{invent-1}
\min\limits_{u_t\geq  0} &\bbe\Big[\, \sum\limits_{t=1}^{T} c_tu_t+\psi_t(x_t+u_t,D_t)\Big] \\
\label{invent-2}
{\rm s.t.} &  x_{t+1}= x_t+u_t-D_t,\;t=1,...,T,
\end{eqnarray}
where $c_t,b_t,h_t$ are the ordering cost, backorder penalty cost and holding cost per unit, respectively, $x_t$ is the current inventory level, $u_t$ is the order quantity,    $D_t$ is the demand   at time $t=1,...,T$, and
$$\psi_t(y_t,d_t):=b_t[d_t- y_t]_+ + h_t[y_t-d_t]_+.$$
 In that format
  $x_t$ are    state variables, $u_t$ are controls, $D_t$ are random disturbances,  and $\U_t:=\bbr_+$ is the   feasibility set (independent of state variables).
At time $t=1$ the  inventory level $x_1$ and demand value $D_0$ are  known (initial conditions). Decision about $u_t$   at stage   $t=1,...,T,$ should be made before a realization of the demand $D_t$ becomes known. That is, at time $t$ history $D_{[t-1]}$  of the demand process is available, but future realizations $D_t,...,D_T$ are not known.

 By making change of variables $y_t=x_t+u_t$,  problem \eqref{invent-1} -\eqref{invent-2} is often written in the following equivalent form
\begin{eqnarray}
\label{invent-1a}
\min\limits_{y_t\geq  x_t} &\bbe\Big[\, \sum\limits_{t=1}^{T} c_t(y_t-x_t)+\psi_t(y_t,D_t)\Big] \\
\label{invent-2a}
{\rm s.t.} &  x_{t+1}= y_t-D_t,\;t=1,...,T.
\end{eqnarray}
Problem \eqref{invent-1a} -\eqref{invent-2a} can be viewed as an SP problem in decision variables $y_t$ and $x_t$.

It is convenient to write   dynamic programming equations for the inventory problem  in the form \eqref{invent-1a} -\eqref{invent-2a}, of course this can be reformulated for the form  \eqref{invent-1} -\eqref{invent-2} as well.
As before we assume here that distribution of the demand process does not depend on our decisions.
At stage $t=T,...,1,$ the  value  functions  satisfy the following equation
\begin{equation}
\label{invent-3}
V_t(x_t, D_{[t-1]})= \inf\limits_{y_{t}\geq x_t} \left\{   c_{t}(y_t-x_t) +
\bbe_{|D_{[t-1]}} \big [\psi_t(y_t,D_t) +V_{t+1}\left (y_t-D_{t}, D_{[t]}\right)  \big ]\right\},
\end{equation}
with $V_{T+1}(\cdot)\equiv 0$. The objective functions $c_t(y_t-x_t)+\psi_t(y_t,D_t)$ are convex in $(x_t,y_t)$ and hence
value functions $V_t(x_t, D_{[t-1]})$  are convex in $x_t$ (see Remark \ref{rem-convsoc}).

The optimal policy  of how much to order at stage $t$ is given by $\bar{u}_t=\bar{y}_t-\bar{x}_t$,   where
\begin{equation}\label{invent-4}
 \bar{y}_t\in \argmin\limits_{y_{t}\geq \bar{x}_t} \left\{   c_{t}(y_t-\bar{x}_t) +
\bbe_{|D_{[t-1]}} \big [\psi_t(y_t,D_t) +V_{t+1}\left (y_t-D_{t}, D_{[t]}\right)  \big ]\right\}.
\end{equation}
That is, starting with the initial value $\bar{x}_1=x_1$,    value $\bar{y}_1$ is computed as an optimal solution of the minimization problem in the right hand side of \eqref{invent-4} for $t=1$   (at the firs stage,  since $D_0$ is deterministic, the conditional expectation in \eqref{invent-4} is just the unconditional expectation with respect to  $D_1$). Consequently  the next inventory level  $\bar{x}_2=\bar{y}_1-D_1$ is computed given an  observed realization of $D_1$.  Then $\bar{y}_2$ is computed as an optimal solution of the minimization problem in the right hand side of \eqref{invent-4} for $t=2$, the inventory level is updated as $\bar{x}_3=\bar{y}_2-D_2$,
and so forth going forward in time. This defines an optimal policy with the inventory level and order quantity  at stage $t$ being functions of observed realization of the history $D_{[t-1]}$ of the demand process.

If  the demand process $D_t$ is stagewise independent, then value functions  $V_t(x_t)$ do not depend on the  history of the demand process, and
equations \eqref{invent-3} become
\begin{equation}
\label{invt-ind}
V_t(x_t)= \inf\limits_{y_{t}\geq x_t} \left\{   c_{t}(y_t-x_t) +
\bbe  \big [\psi_t(y_t,D_t) +V_{t+1}\left (y_t-D_{t}\right)  \big ]\right\},
\end{equation}
with  the                                                                                                                                   expectations taken with  respect to the   marginal distributions of the demand process.
By using convexity of the objective function it is possible to show that in the stagewise independence case  the so-called basestock policy is optimal. That is, optimal control (optimal order quantity)  $\bar{u}_t$ is a function of state $x_t$, and is given by
 \begin{equation}\label{basestock}
\bar{u}_t=\max\{x_t,y^*_t\}-x_t=[y^*_t-x_t]_+,
 \end{equation}
where   $y^*_t$  is the  (unconstrained) minimizer
\begin{equation}\label{minimizer}
 y^*_t\in \argmin\limits_{y_t} \left\{   c_{t}y_t  +
\bbe \big [\psi_t(y_t,D_t) +V_{t+1}\left (y_t-D_{t}\right)  \big ]\right\}.
\end{equation}
That is, if at stage $t$  the current inventory level $x_t$ is greater than or equal to   $y^*_t$, then order nothing; otherwise order $y^*_t-x_t$.
Of course,   computation of  the (critical) values $y^*_t$  is based  on availability of   value functions $V_t(\cdot)$.
}
$\hfill \square$
\end{example}

\setcounter{equation}{0}
\section{Risk Averse and Distributionally Robust   Optimization}
\label{sec-risk}

In the risk averse approach the expectation functional is replaced by a risk measure.
Let $(\O,\F,\bbp)$ be a probability space and $\Z$ be a linear space of measurable functions (random variables)  $Z:\O\to \bbr$. Risk measure $\R$ is a functional assigning a number to random variable $Z\in \Z$, that is $\R:\Z\to\bbr$. In the risk neutral case, $\R$ is the expectation operator, i.e., $\R(Z):=\bbe[Z]=\int_\O Z(\w) d\bbp(\w)$.
In order for this expectation to be well defined we have to restrict the space of considered random variables; for the expectation  it is natural to take the space of integrable functions,  i.e., $\Z:=L_1 (\O,\F,\bbp)$.

It was suggested in  Artzner et al    \cite{ADEH:1999} that reasonable  risk measures should satisfy the following axioms: subadditivity,
  monotonicity,  translation equivariance,  positive homogeneity.    Risk
  measures satisfying these axioms were called {\em coherent}. For a thorough discussion of coherent risk measures we can refer to \cite{follm},\cite{SDR}.
  It is possible to show that if the space $\Z$  is a Banach lattice, in particular if $\Z=L_p(\O,\F,\bbp)$, then
  a (real valued) coherent risk measure  $\R:\Z\to\bbr$  is continuous in the norm topology of $\Z$      (cf., \cite[Proposition 3.1]{RuSh:2006a}). It follows then by the Fenchel - Moreau theorem that $\R$ has the following dual representation
\begin{equation}\label{fmtheorem}
  \R(Z)=\sup_{\zeta\in \cA}\int_\O Z(\w) \zeta(\w) d\bbp(\w), \;\;Z\in \Z,
\end{equation}
where $\cA$ is a convex weakly$^*$  compact  subset of  $\Z^*$. Moreover, the axioms of coherent risk measures imply that   every $\zeta\in \cA$ is a density function, i.e.,    almost surely $\zeta(\w)\ge 0$ and $\int_\O \zeta (\w) d\bbp(\w)=1$  (cf., \cite[Theorem 2.2]{RuSh:2006a}).

An important example of coherent risk measure is the Average Value-at-Risk:
 \begin{eqnarray}
 \avr_\alpha (Z)&:=& (1-\alpha)^{-1}\int_{\alpha}^1 F^{-1}_Z(\tau)d\tau \\
 \label{avr-1}
 &=&\inf_{\theta\in \bbr}\left\{\theta+(1-\alpha)^{-1}\bbe [Z-\theta]_+\right\},\;\;\alpha\in (0,1]
  \label{avr-2}
 \end{eqnarray}
where  $F_Z(z):=\bbp(Z\le z)$ is the  cumulative distribution function (cdf) of $Z$ and
$F^{-1}_Z(\tau):=\inf\{z:F_Z(z)\ge \tau\}$ is the respective  quantile. It is natural in this example to use the space
 $\Z=L_1(\O,\F,\bbp)$ and its dual $\Z^*=L_\infty(\O,\F,\bbp)$.
 In various equivalent forms this risk measure was introduced  in different contexts by different authors under different names, such as Expected Shortfall, Expected Tail Loss, Conditional Value-at-Risk;
    variational form \eqref{avr-2}  was suggested in
 \cite{pflug2000},\cite{ury2}.
 In the dual form, it   has representation \eqref{fmtheorem} with
\[
\begin{array}{l}
\cA=\left \{\zeta\in \Z^*: 0\le \zeta\le (1-\alpha)^{-1},\; \int_\O \zeta d\bbp=1\right\}.
\end{array}
\]

Representation \eqref{fmtheorem} can be viewed from the distributionally robust point of view.
Recall that if $P$ is a probability measure on $(\O,\F)$  which is absolutely  continuous with respect to $\bbp$, then by the Radon - Nikodym Theorem it has density $\zeta=dP/d\bbp$.
 Consider the set $\cM$ of probability measures $P$, absolutely continuous with respect to  $\bbp$, defined as
 $
 \cM:=\{P:dP/d\bbp\in \cA\}.
 $
 Then  representation \eqref{fmtheorem} of  functional $\R$ can be written as\footnote{By writing $\bbe_P[Z]$ we emphasize that the expectation is taken with respect to the probability measure (probability distribution) $P$.}
  \begin{equation}\label{droset-2}
 \R(Z)=\sup_{P\in \cM} \bbe_P[Z].
 \end{equation}

In the Distributionally Robust Optimization (DRO) it  is argued that the ``true" distribution is not known, and consequently a set
$\cM$   of probability   measures (distributions)     is constructed in some way, and referred to as the {\em ambiguity set}. In that framework we can view \eqref{droset-2} as the {\em definition} of the corresponding functional. It is straightforward to verify that the obtained  functional $\R$, defined on an appropriate space of random variables,   satisfies the axioms of coherent risk measures. So we have a certain duality relation between the risk averse and distributionally robust approaches to stochastic optimization.
However, there is an  essential difference   in a way how the corresponding ambiguity set $\cM$ is constructed. In the risk averse approach the set $\cM$ consists of probability measures {\em absolutely continuous} with respect a specified (reference) measure $\bbp$. On the other hand,   in the DRO
 such dominating measure may not be naturally defined. This happens, for instance, in popular settings where  the ambiguity set is defined by moment constraints or is the set of probability measures with prescribed Wasserstein distance  from  the empirical distribution defined by the data.

In order to extend the risk averse and distributionally robust optimization  to the dynamical setting of multistage stochastic optimization, we need to construct a conditional counterpart of the respective  functional $\R$.  In writing dynamic equations \eqref{dynpr-2} and \eqref{soc-2} we used conditional expectations  in a somewhat informal manner.
For the expectation operator there is a  classical definition of conditional expectation. That is,
let $P$ be a probability measure and  $\G$ be  a sigma subalgebra of $\F$. It is said that random variable, denoted $\bbe_{|\G} [Z]$, is the conditional expected value of random variable $Z$ given $\G$ if the following two properties hold: (i) $\bbe_{|\G} [Z]:\O\to\bbr$ is $\G$-measurable, (ii) $\int_A \bbe_{|\G} [Z] (\w)dP(\w)=\int_A Z(\w) dP(\w)$ for any $A\in \G$. There are many versions of $\bbe_{|\G} [Z](\w)$ which can differ from each other on sets of $P$-measure zero.  The conditional expectation depends on measure $P$, therefore we sometimes write $\bbe_{P|\G} [Z]$ to emphasize this when various probability measures are considered.

It seems that it is natural to define the conditional counterpart of the distributionally robust functional, defined in \eqref{droset-2}, as
$\sup_{P\in \cM} \bbe_{P|\G}[Z]$.  However there are technical issues with a precise meaning of such definition since there are different versions of conditional expectations related to different probability measures $P\in \cM$. Even more important, there are  conceptual problems with such definition, and in fact this is not how conditional risk measures are defined (cf., \cite{pichsha}). In the risk averse setting there is a natural concept of law invariance. For a random variable $Z$ its cumulative distribution function (cdf), with respect to the reference probability  measure $\bbp$, is $F_Z(z):=\bbp(Z\le z)$,  $z\in \bbr$.
It is said that random variables $Z$ and $Z'$ are {\em distributionally equivalent} if their
 cumulative distribution functions $F_Z(z)$ and $F_{Z'}(z)$ do coincide for all $z\in \bbr$. It is said that risk measure $\R$ is {\em law invariant} if $\R(Z)=\R(Z')$ for any
distributionally equivalent $Z$ and $Z'$. For example the Average Value-at-Risk measure $\R=\avr_\alpha$ is law invariant.

Note that law invariance is defined with respect to the reference measure $\bbp$. Law invariant coherent risk measure $\R(Z)$ can be considered as a function of  the  cumulative distribution function $F_Z$, and consequently  its conditional counterpart can be defined as the respective function of
the conditional counterpart of  $F_Z$;    this definition  can be made   precise. On the other hand, in the DRO setting where there is no dominating probability measure and the concept of law invariance is not applicable, it is not completely clear how to define a conditional counterpart of functional $\R$ defined in \eqref{droset-2}, although an attempt was made in  \cite{pichsha}. See Remark \ref{rem-risk} below  for a further discussion of law invariance.

In the risk neutral setting we were particularly interested in the stagewise independent case (Assumption \ref{assum-stage}). Under that assumption the dynamic programming equations simplify to \eqref{dynind-1} - \eqref{dynind-1a}  in  the SP and to \eqref{soc-ind} in the  SOC frameworks, respectively. Natural counterparts of these equations are obtained by replacing the expectation operator with a coherent risk measure. That is, in the risk averse  SP framework the counterpart of \eqref{dynind-1a} is
\begin{equation}\label{drisk-1}
  \Q_{t+1}(x_{t}):= \R_{t+1}  [Q_{t+1}(x_t,\xi_{t+1})],\;t=1,...,T-1,
\end{equation}
where  $\R_{t+1}$ is a    coherent risk measure applied to the  random variable   $Q_{t+1}(x_t,\xi_{t+1})$.
For example,
\begin{equation}\label{convcomb}
\R_t(\cdot):=(1-\lambda_t ) \bbe[\,\cdot\,]+\lambda_t \avr_{\alpha_t}(\cdot),\;\;\lambda_t\in [0,1],
\end{equation}
is a law invariant coherent risk measure representing convex combination of the expectation and Average Value-at-Risk.
Such choice of risk measure tries to reach a compromise between minimization of the cost  on average and controlling the risk of upper deviations from the average (upper quantiles  of the cost) at every stage of the process. In applications the parameters $\lambda_t$ and $\alpha_t$ often are constants, independent of the stage $t$, i.e.,  the risk measure $\R_t=\R$ is the same for every stage (cf., \cite{Sha2012a}).

 \begin{remark}
 {\rm 
The  dynamic programming equations,  with \eqref{drisk-1} replacing \eqref{dynind-1a}, correspond
  to the {\em nested} formulation of the   risk averse multistage SP problem (cf., \cite{RuSh:2006b},\cite[Section 6.5]{SDR}).
   In the nested  risk averse multistage optimization the risk is  controlled at {\em every stage} of the process and the total nested value of the cost does not have practical interpretation.  Nevertheless, when solving the corresponding risk averse optimization problem it would be useful to estimate its  optimal value for the purpose of evaluating error of the computed solution.
   } $\hfill \square$
  \end{remark}
 
  In the DRO setting
the  counterpart of equations \eqref{dynind-1a} can be written as
\begin{equation}\label{drisk-2}
  \Q_{t+1}(x_{t}):=  \sup_{P\in \cM_{t+1}}\bbe_P [Q_{t+1}(x_t,\xi_{t+1})],
\end{equation}
where $\cM_{t+1}$ is a set of probability distributions of vector $\xi_{t+1}$.
Note that by convexity and monotonicity properties of coherent risk measures,   convexity of the cost-to-go functions is preserved here. That is, convexity of
  $Q_{t+1}(\cdot,\xi_{t+1})$ implies convexity of $\Q_{t+1}(\cdot)$ in the left hand side of equations \eqref{drisk-1} and \eqref{drisk-2}. This in turn implies convexity of $Q_t(\cdot,\xi_{t})$ in the left hand side of equation \eqref{dynind-1}, provided the function  $f_t(\cdot,\xi_t)$ is convex and the set $\X_t$ is convex. Therefore we have here similar to the risk neutral case (compare with Proposition  \ref{pr-spconvex}) the following convexity property.
\begin{proposition}
\label{pr-con1}
  If the objective
 functions  $f_t(x_t,\xi_t)$ are convex in $x_t$  and the sets $\X_t$ are convex for $t=1,...,T$,
 then the cost-to-go functions $Q_{t+1}(x_t,\xi_{t+1})$ and $\Q_{t+1}(x_t)$ are convex in $x_t$.  In particular such
convexity of the cost-to-go functions follows for risk averse  linear multistage stochastic programs.
\end{proposition}

Similar considerations apply to the  SOC framework, with  the risk averse  counterpart of  equations  \eqref{soc-ind}   written as
\begin{equation}
\label{drisk-3}
V_t(x_t)=\inf\limits_{u_t\in \U_t(x_t)}
\R_t \left [c_t(x_t,u_t,\xi_t)+
V_{t+1}\big(F_t(x_t,u_t,\xi_t) \big)\right],
\end{equation}
where $\R_t$, $t=1,...,T$,  are    coherent risk measures.
Again these dynamic equations correspond to the {\em nested} formulation of the respective SOC problem. Similar to the risk neutral case, the value functions
$V_t(x_t)$ are convex if Assumption \ref{ass-conv} holds.

 \setcounter{equation}{0}
\section{Dynamic Cutting Planes   Algorithms}
\label{sec-cut}

The basic  idea of cutting planes algorithms for solving convex problems is approximation of the (convex) objective function by cutting planes. That is,  let $f:\bbr^n\to \bbr\cup\{+\infty\}$ be an extended real valued function.
Its  domain is $\dom(f):=\{x\in \bbr^n: f(x)<+\infty\}$.
We say that an affine function $\ell(x)=\alpha +\beta^\top x$ is a {\em cutting plane} of $f(x)$ if
$f(x)\ge \ell (x)$ for all $x\in \bbr^n$. If moreover $f(\bar{x})=\ell (\bar{x})$ for some $\bar{x}\in \dom(f)$, then $\ell(x)$ is said to be a {\em supporting plane} of $f(x)$,    at the point $\bar{x}$. When   the function $f$ is convex, its supporting plane at a point $\bar{x}\in \dom(f)$ is given by $\ell(x)=f(\bar{x}) +g^\top (x-\bar{x})$, where $g$ is a subgradient of $f$ at $\bar{x}$. The set of all subgradients, at a point
$x\in \dom(f)$,
 is called the subdifferential of $f$ at $x$  and denoted $\partial f(x)$. The subdifferential of a convex function    at every interior point  of its domain   is  nonempty, \cite[Theorem 23.4]{roc1970}.

For multistage linear SP problems a cutting planes type algorithm, called Stochastic Dual Dynamic Programming (SDDP),  was introduced in Pereira  and Pinto \cite{per1991},  building on nested decomposition algorithm of Birge  \cite{bir85}. The  SDDP  algorithm is going backward and forward in order to  build  piecewise linear  approximations of the   cost-to-go functions by   their cutting planes.

 \begin{remark}
 \label{rem-linear}
 {\rm
 Let us recall how to compute a subgradient of the optimal value function of a linear program. That is, let $Q(x)$ be the optimal value of linear program
  \begin{equation}\label{linpr-1}
 \min_{y\ge 0} c^\top y\;\;{\rm s.t.}\;\;Bx+Ay=b,
 \end{equation}
 considered as a function\footnote{If for some $x$   problem \eqref{linpr-1} does not have a feasible solution, then $Q(x):=+\infty$.}
    of vector $x$.
 The dual of    \eqref{linpr-1}  is  the linear program
    \begin{equation}\label{linpr-2}
 \max_{\lambda}\lambda^\top (b-Bx)\;\;{\rm s.t.}\;\;A^\top\lambda\le c.
 \end{equation}
 Function  $Q(\cdot)$ is  piecewise linear convex, and its
subgradient  at a point $x$ where $Q(x)$ is finite,   is given by $-B^\top \bar{\lambda}$ with  $\bar{\lambda}$  being   an optimal solution of the dual problem \eqref{linpr-2} (e.g., \cite[Section 2.1.1]{SDR}).
 That is,  the subgradient of the optimal value function of a linear program  is computed by solving the dual of that linear program. A more general setting amendable to linear programming formulation, is  
 when the objective function of \eqref{linpr-1}  is  polyhedral.
  In such cases   the  gradient of the respective optimal value function can be computed by solving the dual of the obtained  linear program. This principle is applied in the SDDP algorithm  to approximation of the cost-to-go functions by their cutting planes in a dynamical setting. This motivates the  ``Dual Dynamic" part  in the name of the algorithm.
 } $\hfill \square$
 \end{remark}

 \setcounter{equation}{0}
\subsection{SDDP algorithm for SP problems }
\label{sec-stprcut}

Consider SP problem \eqref{stpr-1} - \eqref{stpr-1a}. Suppose  that   Assumptions \ref{assum-ind} and \ref{assum-stage} hold, and hence the corresponding cost-to-go functions are defined by equations \eqref{dynind-1} - \eqref{dynind-1a}. If probability distributions of random vectors $\xi_t$ are continuous, these distributions should be  discretized in order to compute the respective expectations (see Remark \ref{rem-cont} below).
So  we also make the following assumption.

\begin{assumption}
\label{ass-finite}
The  probability distribution of $\xi_t$ has  finite support $\{\xi_{t1},...,\xi_{tN}\}$ with respective probabilities\footnote{For the sake of simplicity we assume that the number $N$  of realizations of $\xi_t$  is the same for every stage $t=2,...,T$. Recall that in the SP framework,  $\xi_1$ is deterministic.}  $p_{t1},...,p_{tN}$, $t=2,...,T$.
Denote $A_{tj}:=A_t(\xi_{tj})$, $B_{tj}:=B_t(\xi_{tj})$, $b_{tj}:=b_t(\xi_{tj})$, and
$Q_{tj}(x_{t-1}):=Q_{t}(x_{t-1},\xi_{tj})$  the cost-to-cost-to-go functions,
 defined in \eqref{dynind-1}.
 \end{assumption}

The corresponding  expected value
cost-to-go functions
 can be written as
\begin{equation}
\label{dynex}
\begin{array}{ll}
  \Q_{t}(x_{t-1})=\sum_{j=1}^N p_{tj} Q_{tj}(x_{t-1}).
  \end{array}
\end{equation}
Suppose further that the objective  functions $f_t(x_t,\xi_{tj}):=c_{tj}^\top x_t$  are linear in $x_t$,    and the sets $\X_t:=\bbr^{n_t}_+$, i.e., the constraint $x_t\in \X_t$ can be written as $x_t\ge 0$, $t=1,...,T$. In that case  \eqref{stpr-1} - \eqref{stpr-1a} becomes the standard linear multistage SP problem:
\begin{equation}\label{stplin}
\min\limits_{x_t\ge 0}  \bbe\Big[\sum_{t=1}^T c_t^\top x_t\Big]\;\;{\rm s.t.}\;
B_t x_{t-1} +A_t x_t=b_t,\;t=1,...,T.
\end{equation}

It could be noted that if the sets $\X_t$ are polyhedral and the  functions $f_t(x_t,\xi_t)$ are polyhedral in $x_t$, then the problem is linear and can be formulated  in the standard form by the  well known    techniques of linear programming.

\begin{remark}
\label{rem-cont}
{\rm
From the modeling point of view, random vectors $\xi_t$ often are assumed to have continuous distributions. In such cases these continuous distributions have to be discretized  in order to perform numerical calculations. One possible approach to such discretization is to use Monte Carlo sampling techniques. That is, a random sample $\xi_{t1},...,\xi_{tN}$ of size $N$ is generated from the probability distribution of random vector $\xi_t$, $t=2,...,T$. This approach is often referred to as the Sample Average Approximation (SAA) method.
The constructed discretized problem can be considered in the above  framework with equal probabilities $p_{tj}=1/N$, $j=1,...,N$.
Statistical properties of the SAA method in the multistage setting are discussed, e.g., in \cite[Section 5.8]{SDR}. It  could be mentioned that the question of the sample size $N$, i.e. how many discretization points  per stage to use, can be only addressed in the framework of the original  model with continuous distributions. One of the advantages of the SAA method is that statistical tests can be applied to the results of computational procedures, this is discussed   in \cite{ding2019}.
}
$\hfill \square$
\end{remark}

The SDDP algorithm for linear multistage SP problems, having finite number of scenarios,  was described in details in a number of publications  (see, e.g., recent survey paper \cite{FR2021} and references therein). For a  later reference   we briefly discuss it below.
The SDDP algorithm consists of backward and forward steps. At each iteration of the algorithm,  in the  backward step   current approximations of the cost-to-go functions  $\Q_{t}(\cdot)$, by cutting planes,  are updated by adding  respective cuts. That is, at the last stage the cost-to-go function  $Q_{Tj}(x_{T-1})$, $j=1,...,N$, is given by the optimal solution of the linear program
\begin{equation}\label{sddp-1}
\min _{x_T\ge 0} c_{Tj}^\top x_T\;\;{\rm s.t.}\;\;B_{Tj}x_{T-1}+A_{Tj}x_{T}=b_{Tj}.
\end{equation}
The subgradient   of $Q_{Tj}(x_{T-1})$ at a chosen point $\tilde{x}_{T-1}$ is given by $g_{Tj}=-B_{Tj}\tilde{\lambda}_{Tj}$,   where
$\tilde{\lambda}_{Tj}$   is an optimal solution of the dual of  problem \eqref{sddp-1} (see Remark \ref{rem-linear} above). Consequently the   subgradient of $\Q_{T}(x_{T-1})$ at   $\tilde{x}_{T-1}$ is computed as
$\sg_T=\sum_{j=1}^N p_{Tj}g_{Tj}$, and hence the
cutting plane
\begin{equation}\label{sddp-2}
 \ell_T(x_{T-1})=\Q_{T}(\tilde{x}_{T-1})+\sg_T^\top (x_{T-1}-\tilde{x}_{T-1})
\end{equation}
is added to the current family of cutting  planes of $\Q_{T}(\cdot)$.
Actually at the last stage the above cutting plane $\ell_T(\cdot)$ is the {\em supporting}  plane of $\Q_{T}(\cdot)$ at $\tilde{x}_{T-1}$. This is because
$\Q_{T}(\tilde{x}_{T-1})$ is explicitly computed by solving problem \eqref{sddp-1}.

At one step back at stage $T-1$, the   cost-to-go function  $Q_{T-1,j}(x_{T-2})$, $j=1,...,N$, is given by the optimal solution of the   program
\begin{equation}\label{sddp-3}
\min _{x_{T-1}\ge 0} c_{T-1,j}^\top x_{T-1}+\Q_T(x_{T-1})\;\;{\rm s.t.}\;\;B_{T-1,j}x_{T-2}+A_{T-1j}x_{T-1}=b_{T-1,j}.
\end{equation}
The cost-to-go function  $\Q_T(\cdot)$ is not known. Therefore it is  replaced by the  approximation $\und_T(\cdot)$ given by the maximum of the current  family of cutting planes including the cutting plane added at stage $T$ as described above. Since $\und_T(\cdot)$ is the maximum of a  finite number of affine functions,  the obtained  problems
\begin{equation}\label{sddp-4}
\min _{x_{T-1}\ge 0} c_{T-1,j}^\top x_{T-1}+\und_T(x_{T-1})\;\;{\rm s.t.}\;\;B_{T-1,j}x_{T-2}+A_{T-1j}x_{T-1}=b_{T-1,j},
\end{equation}
$j=1,...,N$, are  linear. Let $\cQ_{T-1,j}(x_{T-2})$ be the optimal value function of problem
\eqref{sddp-4} and $\bar{\cQ}_{T-1}(x_{T-2})=\sum_{j=1}^N p_{T-1,j}\cQ_{T-1,j}(x_{T-2})$ be the corresponding  estimate of the   expected value cost function. Note that since $\Q_{T}(\cdot)\ge \und_T(\cdot)$, we have that  $Q_{T-1,j}(\cdot) \ge \cQ_{T-1,j}(\cdot)$, $j=1,...,N$, and hence  $\Q_{T-1}(\cdot) \ge \bar{\cQ}_{T-1}(\cdot)$.
Gradient $g_{T-1,j}$ of the optimal value function $\cQ_{T-1,j}(x_{T-2})$    can be computed at a chosen point $\tilde{x}_{T-2}$ by solving the dual of problem \eqref{sddp-4} and hence computing its optimal solution. Then
  $\sg_{T-1}=\sum_{j=1}^N p_{T-1,j}g_{T-1,j}$ is the gradient of  $\bar{\cQ}_{T-1}(x_{T-2})$ at $\tilde{x}_{T-2}$. Consequently
\begin{equation}\label{sddp-5}
 \ell_{T-1}(x_{T-2})=\und_{T-1}(\tilde{x}_{T-2})+\sg_{T-1}^\top (x_{T-2}-\tilde{x}_{T-2})
\end{equation}
is a  cutting plane of $\Q_{T-1}(\cdot)$, and
is added to the current collection of cutting planes of $\Q_{T-1}(\cdot)$.
And so on going backward in time until the  corresponding  approximation of the first stage problem is obtained:
\begin{equation}\label{sddp-6}
\min _{x_1\ge 0} c_1^\top x_{1}+\und_2(x_{1})\;\;{\rm s.t.}\;\; A_{1}x_{1}=b_{1}.
\end{equation}
Let us note that by the construction, $\Q_t(\cdot)\ge \und_t(\cdot)$ for $t=2,...,T$. Therefore any cutting plane of
$\und_t(\cdot)$ is also a cutting plane of $\Q_t(\cdot)$.
Also it follows  that the optimal value of problem \eqref{sddp-6} is less than or equal to the optimal value of  problem \eqref{stpr-1} - \eqref{stpr-1a}. Therefore   optimal value of problem \eqref{sddp-6}
provides a (deterministic) lower bound for the optimal value  the considered multistage problem.

The computed estimates $\und_t(\cdot)$, $t=2,...,T$, of the cost-to-go functions, and optimal solution $\bar{x}_1$ of the   first stage problem \eqref{sddp-6},
define a feasible policy in the same way as in \eqref{dynind-2}. Since this policy is feasible, its value  is greater than or equal to the optimal value of the considered problem \eqref{stpr-1} - \eqref{stpr-1a}. In the forward step this value is estimated by randomly generating  realizations of the data process and evaluating the corresponding policy values. That is, consider
a  sample path (scenario) $\hat{\xi}_2,...,\hat{\xi}_T$ of the data process.
Let $\hat{B}_t=B_t(\hat{\xi}_t)$, $\hat{A}_t=A_t(\hat{\xi}_t)$, $\hat{b}_t=b_t(\hat{\xi}_t)$ and $\hat{c}_t=c_t(\hat{\xi}_t)$ be  realizations of the random parameters corresponding to the generated  sample path.
Starting  with first stage solution  $\bar{x}_1$, next values $\bar{x}_t$ are computed iteratively going forward in time by computing an optimal solution of the optimization problem in  \eqref{dynind-2} with $\Q_{t+1}(\cdot)$ replaced by $\und_{t+1}(\cdot)$. That is
  \begin{equation}\label{polappl}
 \bar{x}_t \in   \argmin_{x_t\ge 0}   \left\{ \hat{c}_t^\top x_t +\und_{t+1}(x_{t}): \hat{B}_t \bar{x}_{t-1} +\hat{A}_t x_t=\hat{b}_t\right\}.
\end{equation}
  The corresponding  total sum
  \begin{equation}\label{value}
    \hat{\vv}:=\hat{c}_1^\top \bar{x}_1+\sum_{t=2}^T  \hat{c}_t^\top \bar{x}_t
  \end{equation}
  is a function of the generated sample   and hence is random. Its expected value $\bbe[\hat{\vv}]$
  is equal to the value of the constructed   policy, and thus
  $\bbe[\hat{\vv}]$ is greater than or equal to the optimal value of   problem \eqref{stpr-1} - \eqref{stpr-1a}.

   The forward step of the algorithm proceeds to evaluate value of the constructed policy by randomly generating the sample paths (scenarios). At every stage $t=2,...,T$, a point $\xi_t$ is generated at random according to the distribution of the corresponding random vector.  That is,  a point  $\hat{\xi}_t$ is chosen with probability $p_{ti}$, $i=1,...,N$,
from the set $\{\xi_{t1},...,\xi_{tN}\}$.
This generates the corresponding   sample path (scenario)
$\hat{\xi}_2,...,\hat{\xi}_T$ of the data process. Consequently the policy values for the generated sample path are computed, by using the rule \eqref{polappl},
and thus  the   total sum $\hat{\vv}$, defined in \eqref{value}, is evaluated. The computed value $\hat{\vv}$
gives a point estimate of   value of the  constructed policy.
This procedure can be repeated several times   by randomly (and independently from each other) generating   $M$   scenarios. Let $\hat{\vv}^1,...,\hat{\vv}^M$ be values,  calculated in accordance with \eqref{value}, for each generated scenario. Their average
$\bar{\vv}:=M^{-1}\sum_{i=1}^M \hat{\vv}^i$ gives an unbiased estimate of the value of the considered policy. An upper edge  $\bar{\vv}+z_\alpha S/\sqrt{M}$  of the corresponding confidence interval is viewed as a statistical estimate of an upper bound for the optimal value of the considered multistage problem. Here  $S^2=(M-1)^{-1}\sum_{i=1}^M (\hat{\vv}^i-\bar{\vv})^2$ is an   estimate of the variance $\var (\hat{\vv})$, and $z_\alpha$ is the critical level, say $z_\alpha=2$. The justification for this bound is   that by the Central Limit Theorem (CLT) for large $M$, the average $\bar{\vv}$ has approximately normal distribution  and $z_\alpha$ corresponds to mass $100(1-\alpha)\%$ of the standard normal distribution. For $M$ not too large it is a common practice to use  critical values from $t$-distribution with $M-1$ degrees of freedom, rather than normal.  In any case    this is just a heuristic justified by  the CLT. Nevertheless it is useful and typically gives a reasonable estimate of the upper bound.

\begin{remark}
{\rm
The constructed approximations $\und_t(\cdot)$, of the expected value  cost-to-go functions, also define a feasible policy for the original problem with continuous distributions of the data process  vectors $\xi_t$.  Its value can be estimated in the same way by generating a sample path  $\hat{\xi}_2,...,\hat{\xi}_T$ from the {\em original}  continuous distribution and   computing the corresponding  policy  values in the forward step in accordance with \eqref{polappl}. Consequently the policy value is estimated by averaging the computed values for a number of randomly generated sample paths. This can be used for construction of statistical upper bound for the original problem. Expectation of the optimal value   of the SAA problem, based on randomly generated sample, is less than or equal to the optimal value of the original problem. Therefore the lower bound generated by the SDDP algorithm  gives a point estimate of a lower bound of the original problem.
We can refer to \cite{ding2019} for various related  statistical testing procedures.
} $\hfill \square$
\end{remark}

The difference between the (statistical) upper bound and (deterministic) lower bound provides an upper bound for  the optimality gap of the constructed policy. It can be used as a stopping rule by stopping the iteration procedure when this estimated gap is smaller than a specified precision level. The bound is quite conservative and for larger problems may not be reducible  to the specified accuracy  level in a reasonable computational time. A more practical heuristic approach is to stop the iterations after the lower bound stabilizes and new iterations, with additional cutting planes, do not significantly improve (increase) the lower bound.

 The above procedure depends on a choice of the  points $\tilde{x}_t$, referred to as the  {\em trial points}.
The  forward step   suggests to use the computed values $\bar{x}_t$ as trial points  in the next iteration of the backward step of the algorithm. The rational for this can be explained as follows. In general,  in order to approximate the cost-to-go functions {\em uniformly}  with a given accuracy in a specified region, the number of cutting planes grows exponentially with increase of the dimension of decision variables $x_t$. In many applications the regions where the cost-to-go functions should be approximated more accurately are considerably smaller than  the   regions of possible variability of decision variables.
By randomly generating scenarios,  the forward step   generates trial points where likely the cost-to-go functions should be approximated more accurately.

\begin{remark}
\label{rem-regul}
{\rm
  The trial points are supposed to be chosen at stages $t=2,...,T-1$.
This suggests that the above procedure is relevant for number of stages $T\ge 3$.  For two stage problems, i.e., when $T=2$, there is no sampling involved  and  this becomes the classical Kelley's cutting plane algorithm, \cite{kel:60}.
 It is well known in the deterministic convex optimization that Kelley's algorithm can behave poorly with increase of the number of decision variables, and much better cutting plane algorithms can be designed by  using regularization or trust region/level set strategies~\cite{LNN1995,LAN2015}. However, it is not clear how such methods can be extended to the dynamic multistage settings (although such attempts were made). The reason is that in the multistage problems the optimization is performed over policies which are functions of the data process. It is not known  a priori  where   realizations of considered policies may happen in future stages and an attempt to restrict  regions of search by some form of regularization could be not  helpful.
} $\hfill \square$
\end{remark}

Proofs of convergence of the SDDP type algorithms were studied in several publications (e.g., \cite{Gir:2014} and references there in).  These  proofs are based on the argument that eventually almost surely   the forward step of the algorithm will go through every possible scenario many times, and hence the cost-to-go functions will be approximated with any accuracy at relevant regions. Since   number of scenarios grows exponentially with increase of the number of stages, this is not a very realistic assumption.  Numerical experiments indicate that,  in a certain sense,  the numerical complexity of SDDP method grows almost linearly with  increase of the number of stages while the dimension of decision vectors is constant.  This  property, of the rate of convergence, was investigated in recent paper \cite{Lan22}.

\subsubsection{Interstage dependence}
\label{sec-intdep}

Suppose now that the data process $\xi_t\in \bbr^d$ is Markovian, and hence the cost-to-go functions are defined by equations \eqref{dynind-3} - \eqref{dynind-3a}.
There are different ways how such Markovian structure can be modeled. One classical approach is time series analysis. The data process is modeled as, say first order,  autoregressive process, i.e.,
\begin{equation}\label{tsaut-1}
  \xi_{t}=\mu+\Phi \xi_{t-1} +\e_t,
\end{equation}
where $\Phi\in \bbr^{d\times d}$ and $\mu\in \bbr^d$ are parameters of the process and $\e_t$ is a sequence of random error vectors. It is assumed that
the errors process $\e_t$ is   stagewise independent (in fact it is usually assumed that $\e_t$ are  i.i.d.). Suppose further that  only the right hand side variables  $b_t$  in the considered problem
 are  random and set $\xi_t:=b_t$.
Then we can write the cost-to-go functions in terms of $x_t$ and $\xi_t$ treating $\e_t$ as the corresponding random process:
\begin{equation}\label{tsaut-2}
Q_t(x_{t-1},\xi_{t-1},\e_{t})=  \inf_{x_t\ge 0,\,\xi_t} \left\{ c_t^\top x_t + \Q_{t+1}(x_{t},\xi_{t}):
\begin{array}{lll}
B_t  x_{t-1} +A_t x_t=\xi_t,\\
\xi_{t}-\mu-\Phi \xi_{t-1} =\e_t
\end{array}
\right\},
\end{equation}
with
\begin{equation}\label{tsaut-3}
\Q_{t+1}(x_{t},\xi_{t})=  \bbe [Q_{t+1}(x_t,\xi_{t},\e_{t+1})].
\end{equation}
Here $(x_{t},\xi_{t})$ are treated as decision variables, $\e_t$ is viewed as the corresponding random process with the expectation in \eqref{tsaut-3} is taken with respect to $\e_{t+1}$.

There are two drawback in the above approach. First, the number of decision variables is artificially increased with the cutting planes in the corresponding SDDP algorithm computed jointly in $x_t$ and $\xi_t$. This increases the computational complexity of the numerical procedure. Second, in order  to have a {\em linear} problem, only the right hand side parameters are allowed to be random. An alternative approach is to approximate  Markovian process $\xi_t$ by Markov chain, which we are going to discuss next.

Suppose that distribution of random vector $\xi_t$ is supported on set  $\Xi_t\subset \bbr^{d_t}$, $t=2,...,T$. Consider the partition\footnote{For the sake of simplicity we assume that the number $N$ of cells is the same at every stage.}
$\Xi_{ti}$,  $i=1,..., N$,  of $\Xi_t$ by Voronoi cells, i.e.,
$$\Xi_{ti}:=\left\{\xi_t\in \Xi_t:\|\xi_t-\eta_{ti}\|\le \|\xi_t-\eta_{tj}\|\;\text{for all}\;j\ne i\right\},$$
where $\eta_{ti}\in \Xi_t$ are given points. That is, $\Xi_{ti}$ is the set of points of $\Xi_t$ closest to $\eta_{ti}$. There are various ways how the center points $\eta_{ti}$  can be chosen. For example,
   $\eta_{ti}$  could be determined by
minimizing  the squared distances:
\begin{equation}\label{means}
 \min_{\eta_{t1},...,\eta_{t N}}\int\min_{i=1,...,N}\|\xi_t-\eta_{ti}\|^2dP(\xi_t).
\end{equation}
 The above optimization  problem is not convex. Nevertheless it can be reasonably solved by various methods (cf., \cite{Lohndorf}).
Consider conditional probabilities $p_{tij}:=P(\xi_{t+1}\in \Xi_{t+1,j}|\xi_t\in \Xi_{ti})$. Note that different Voronoi cells can have  common boundary points, we assume that probabilities of such events are zero so that the probabilities $p_{tij}$ are well defined.
 This leads to an approximation of the data process by a (possibly nonhomogeneous)  Markov chain with transition probabilities $p_{tij}$ of moving from point $\eta_{ti}$ to point $\eta_{t+1,j}$ at the next stage.

  Consider $c_{ti}:=c_{t}(\eta_{ti})$, $A_{ti}:=A_{t}(\eta_{ti})$, $B_{ti}:=B_{t}(\eta_{ti})$ and $b_{ti}:=b_{t}(\eta_{ti})$,  $i=1,...,N$. Recall dynamic equations \eqref{dynind-3} - \eqref{dynind-3a}  for   Markovian data process.
 For the constructed Markov chain these dynamic equations   can be written as
\begin{equation}\label{tsaut-4}
Q_t(x_{t-1},\eta_{ti})=  \inf_{x_t\ge 0} \Big\{ c_{ti}^\top x_t + \Q_{t+1}(x_t,\eta_{ti}): B_{ti}  x_{t-1} +A_{ti} x_t=b_{ti}\Big\},\;i=1,...,N,
\end{equation}
  where
\begin{equation}\label{tsaut-5}
\Q_{t+1}(x_t,\eta_{ti})=\sum_{j=1}^N p_{tij} Q_{t+1}(x_t,\eta_{t+1,j}).
\end{equation}
Note that  the cost-to-go functions  $Q_t(x_{t-1},\eta_{t})$  are convex in
$x_{t-1}$, while the argument  $\eta_{t}$ can take only the values
$\eta_{t1},...,\eta_{tN}$.

It is possible to apply the SDDP method separately  for every $\eta_{t}=\eta_{ti}$, $i=1,...,N$ (cf.,  \cite{Lohndorf}).
That is, in the backward step of the algorithm the cost-to-go functions
$\Q_{t+1,i}(x_t):=\Q_{t+1}(x_t,\eta_{ti})$    are approximated by their cutting planes (in $x_t$) separately for each $\eta_{ti}$, $i=1,...,N$.
This defines a policy for every realization (sample path) $\eta_1,...,$  of the constructed  Markov chain.  This  can be used in the forward step of the algorithm.

  It could happen that a sample path  $\xi_1,...,$  of the original data process is different from every sample  path $\eta_1,...,$ of the constructed Markov chain.
In order   to extend the constructed policy  to a policy for the original data process, one approach is to define the policy values for a sample path of the original data process by taking the corresponding policy values from the sample path     of the Markov chain   closest in some sense  to the considered  sample path  of the original data process.

The above approach is based on discretization (approximation)  of the Markov process by a Markov chain.  Compared with the   approach based on the  autoregressive modeling of the data process, the approach of Markov chain discretization   is not restricted to the case where only the  right hand side parameters are random. Also it uses cutting planes approximations only with respect to variables $x_t$, which could be computationally advantageous. On the other hand, its computational complexity and required computer memory grows with increase of the number $N$ of   per stage discretization points. Also   unlike the SAA method, there is no available  inference describing statistical properties of  such discretization.  From an intuitive point of view,   such discretization could become  problematic with increase of the dimension of the random vectors of the data process.

\subsubsection {Duality   of multistage linear stochastic
programs}
\label{sec-dual}

The linear multistage stochastic program \eqref{stplin}, with a finite number of scenarios, can be viewed as a large scale linear program, the so called deterministic equivalent. By the standard theory of linear programming it has the (Lagrangian) dual.
Suppose that the primal problem \eqref{stplin} has finite optimal value.  Then the optimal values of the primal and its dual are equal to each other and both problems have optimal solutions. Dualization of the feasibility constraints  leads to the following (Lagrangian) dual of problem \eqref{stplin} (cf., \cite[Section 3.2.3]{SDR})
\begin{equation}\label{eq-2}
\begin{array}{cll}
  \max\limits_{\pi} & \bbe\big[  \sum_{t=1}^T     b_t^\top \pi_t\big]\\
    {\rm s.t.}  &A^\top_T\pi_T\le c_T,\\
  &  A_{t -1}^\top \pi_{t-1}+
    \bbe_{|\xi_{[t-1]}}\left[ B_{t}^\top \pi_{t}\right ]\le c_{t-1},\;t=2,...,T.
   \end{array}
\end{equation}
The optimization in \eqref{eq-2} is over policies $\pi_t=\pi_t(\xi_{[t]})$, $t=1,...,T$.

  It is possible to write dynamic programming equations   for the dual problem.
 The idea of applying the SDDP algorithm to the dual problem was introduced in
 \cite{Leclere}. In what follows we use the approach of
  \cite{GSC22} based on formulation \eqref{eq-2}. As before we make Assumptions  \ref{assum-ind},   \ref{assum-stage} and \ref{ass-finite}.
    At the last stage $t=T$,
given $\pi_{T-1}$ and
$\xi_{[T-1]}$, we need to solve the following  problem with respect to $\pi_T$:
\begin{equation}\label{eq-3a-1}
\begin{array}{cll}
  \max\limits_{\pi_{T 1},\ldots,
  \pi_{T N}} &  \sum\limits_{j=1}^{N}   p_{T j}  b_{T j}^\top \pi_{T j}\\
    {\rm s.t.}&
  A_{T j}^\top \pi_{T j}  \le c_{T j},\;j=1,...,N,\\
    &A_{T-1}^\top \pi_{T-1}+
  \sum\limits_{j=1}^{N}   p_{T j}
  B_{T j}^\top \pi_{T j} \le c_{T-1}.
   \end{array}
\end{equation}
Note that  $\sum\limits_{j=1}^{N}   p_{T j}  b_{T j}^\top \pi_{T j}=\bbe[  b_T^\top \pi_T]$ and $\sum\limits_{j=1}^{N}   p_{T j}
  B_{T j}^\top \pi_{T j}  =\bbe \left[ B_{T}^\top \pi_{T}\right ]$.
Since $\xi_T$  is independent of $\xi_{[T-1]}$,  the expectation in \eqref{eq-3a-1} is unconditional with respect to the distribution of $\xi_T$.

The optimal value $V_T(\pi_{T-1},\xi_{T-1})$ and an optimal   solution
$(\bar{\pi}_{T 1},\ldots,\bar{\pi}_{TN})$  of   problem \eqref{eq-3a-1}  are  functions of vectors  $\pi_{T-1}$ and $c_{T-1}=c_{T-1}(\xi_{T-1})$ and matrix  $A_{T-1}=A_{T-1}(\xi_{T-1})$.
And  so on going backward in time, using  the stagewise independence assumption,  we can write the respective dynamic programming equations for
  $t=T-1,...,2$, as
 \begin{equation}\label{eq-6-1}
\begin{array}{cll}
  \max\limits_{\pi_{t 1},\ldots,\pi_{t N}} &     \sum\limits_{j=1}^{N}
  p_{t j}  \left[ b_{t j}^\top \pi_{t j}+ V_{t+1}(\pi_{t j},\xi_{t j})\right]\\
    {\rm s.t.}&
  A_{t-1}^\top \pi_{t-1}+
  \sum\limits_{j=1}^{N} p_{t j} B_{t j}^\top  \pi_{t j}  \le c_{t-1},
   \end{array}
\end{equation}
with $V_t(\pi_{t-1},\xi_{t-1})$ being the optimal value  of problem \eqref{eq-6-1}.
Finally at the first stage the following problem should be solved\footnote{Value functions $V_t$ here should not be confused with the value functions of the SOC framework.}
 \begin{equation}\label{eq-7-1}
  \max_{\pi_{1}}      b_1^\top \pi_1+ V_{2}(\pi_1,\xi_1) .
   \end{equation}

  Let us make the following observations  about   the dual problem.
  Unlike    the primal problem,
  the optimization (maximization) problems \eqref{eq-3a-1} and \eqref{eq-6-1}   do not decompose into separate problems with respect to each $\pi_{t j}$ and should be solved as one linear program with respect to $(\pi_{t 1},...,\pi_{t N})$.
If $A_t$ and $c_t$, $t=2,...,T$,  are deterministic, then $V_t(\pi_{t-1})$ is only a function  of $\pi_{t-1}$.
 The value (cost-to-go) function   $V_t(\pi_{t-1},\xi_{t-1})$   is a concave function of $\pi_{t-1}$. This allows to apply an SDDP type algorithm to the dual problem (see  \cite{GSC22} for details).    For any feasible policy of the dual problem, value of the first stage problem \eqref{eq-7-1} gives an upper bound for the optimal value of the primal problem \eqref{stplin}. Since the dual   is a maximization problem, the SDDP algorithm applied to the dual problem   produces  a deterministic   upper  bound for the optimal value of  problem \eqref{stplin}.

  It is interesting to note that it could happen that the dual problem does  not have relatively complete recourse even if the primal problem has it. It is said that the problem (primal or dual) has   {\em relatively complete recourse},  if at every stage $t=2,...,T$,
for any generated solutions  by the forward process at the previous stages, the respective dynamic program has a feasible
solution at the stage $t$ for every realization of the random data. Without relatively complete recourse it could happen at a certain  stage   that for some realizations of the data process, the corresponding optimization problem is infeasible and the forward process cannot continue. Note that the relatively complete recourse is related to the policy  construction in the forward dynamic process and may not hold even if the problem has finite optimal value. It is still possible to proceed with  construction of the dual upper bound by introducing a penalty term to the dual problem (cf.,  \cite{GSC22}).

\subsection{Cutting planes  algorithm for SOC problems }
\label{sec-socsddp}

Consider the SOC problem \eqref{soc} - \eqref{soc-b}. Suppose that Assumptions \ref{assum-ind} and \ref{assum-stage} hold,  and hence equation \eqref{soc-ind} for the value function follows. Suppose further that Assumption
\ref{ass-conv} is satisfied,  and hence value functions $V_t(\cdot)$ are convex. Then a cutting planes algorithm of the SDDP type can be applied in a way similar to the SP setting. There are, however, some subtle differences which we are going to discuss next.

In the backward step of the algorithm an approximation $\undv_t(\cdot)$ of the value function $V_t(\cdot)$  is constructed by maximum of a family of cutting planes, similar to the SP setting, going backward in time. Let $\tilde{x}_t$, $t=2,...,T+1$,  be trial points, generated say by the forward step of the algorithm at the previous iteration. At the backward step,  starting at  time $t=T+1$  the cutting (supporting)  plane $\ell_{T+1}(x_{T+1})=c_{T+1}(\tilde{x}_{T+1})+\sg_{T+1}^\top (x_{T+1}-\tilde{x}_{T+1})$, where $\sg_{T+1}$ is a subgradient of $c_{T+1}(\cdot)$ at $\tilde{x}_{T+1}$,
is added to the current family of cutting planes of $V_{T+1}(\cdot)$.  Going backward in time,  at stage $t$  we need to compute a subgradient of function
\begin{equation}\label{socsddp-1}
\inf\limits_{u_t\in \U_t(x_t)}
\sum_{j=1}^N p_{tj} \left [c_{tj}(x_t,u_t)+
\undv_{t+1}\big(A_{tj}x_t +B_{tj}u_t+b_{tj}  \big)\right]
\end{equation}
at the trial point $x_t=\tilde{x}_t$, where  $\undv_{t+1}(\cdot)=\max_{1\le i\le M}\ell_{t+1,i}(\cdot)$ is the current approximation of $V_{t+1}(\cdot)$ by cutting planes $\ell_{t+1,i}(\cdot)$. The minimization problem \eqref{socsddp-1} is linear   if functions $c_{tj}(x_t,u_t)$  are linear (polyhedral) in $u_t$,   and $\U_t(x_t)$ is of the form \eqref{setrepr} with the corresponding  mapping $g_t(x_t,u_t)$  being affine in $u_t$.

Let us first consider the case where  $\U_t(x_t)\equiv \U_t$ does not depend on $x_t$. Then the required subgradient is given by (cf., \cite[Theorem 23.8]{roc74})
\begin{equation}\label{subd-3a}
\sg_t=\sum_{j=1}^N  p_{tj}  \left [\nabla  c_{tj}(\tilde{x}_t,\tilde{u}_t)+
A^\top_{tj}  \nabla \undv_{t+1}\big(A_{tj} \tilde{x}_t+B_{tj} \tilde{u}_t+ b_{tj} \big)\right],
\end{equation}
where   $\tilde{u}_t\in \U_t$  is a minimizer of  problem \eqref{socsddp-1} for $x_t=\tilde{x}_t$,  and   $\nabla   c_{tj}(\tilde{x}_t,\tilde{u}_t)$ is a subgradient of $c_{tj}(\cdot,\tilde{u}_t)$ at $\tilde{x}_t$. We also need to
compute a subgradient of $\nabla \undv_{t+1}(\cdot)$ at
the points $A_{tj} \tilde{x}_t+B_{tj} \tilde{u}_t+ b_{tj}$, $j=1,...,N$.
The subgradient  $\nabla \undv_{t+1}(\cdot)$ at a point  $x_{t+1}$  is given by the gradient  $\nabla \ell_{t+1,\nu} (x_{t+1})$ of its cutting plane with    $\nu\in \{1,...,M\}$  such that  $\ell_{t+1,\nu} (x_{t+1})=\undv_{t+1}(x_{t+1})$, i.e.,  $\ell_{t+1,\nu}$ is the supporting plane of  $\undv_{t+1}(\cdot)$ at $x_{t+1}$.
Of course, the gradient of affine function $\ell (x)=\alpha+\beta^\top x$ is $\nabla \ell (x)=\beta$ for any $x$. The cutting plane is then $\underline{V}_t (\tilde{x}_t)+\sg_t^\top(x_t-\tilde{x}_t)$, where
\[
\underline{V}_t (\tilde{x}_t)= \sum_{j=1}^N p_{tj} \left [c_{tj}(\tilde{x}_t,\tilde{u}_t)+
\undv_{t+1}\big(A_{tj}\tilde{x}_t +B_{tj}\tilde{u}_t+b_{tj}  \big)\right].
\]

 It could  be noted that in the present case the required subgradient, and hence the corresponding cutting plane,  can be computed without solving the dual problem. Therefore it may be inappropriate in this setting to call it Stochastic {\em Dual} Dynamic Programming.

When $\U_t(x_t)$ depends on $x_t$, calculation of a subgradient of function \eqref{socsddp-1} could be more involved; nevertheless  in some cases it is still possible to proceed. Suppose that $\U_t(x_t)$ is of the form \eqref{setrepr} with the corresponding  mapping $g_t(\cdot,\cdot)$ having convex components. That is, minimization  problem \eqref{socsddp-1} can be written as
\begin{equation}\label{socmin}
\min\limits_{u_t}
\sum_{j=1}^N p_{tj} \left [c_{tj}(x_t,u_t)+
\undv_{t+1}\big(A_{tj}x_t +B_{tj}u_t+b_{tj}  \big)\right]\;\;{\rm s.t.}\;g_t(x_t,u_t)\le 0.
\end{equation}
Then the required gradient is given by the gradient of the Lagrangian $L_t(x_t,\bar{u}_t,\bar{\lambda}_t)$ of problem \eqref{socmin} with $\bar{u}_t$
being the respective optimal solution and $\bar{\lambda}_t$ being the corresponding Lagrange multipliers vector, provided the optimal solution and Lagrange multipliers are unique (e.g., \cite[Theorem 4.24]{BS2000}).

The forward step of the algorithm is performed similarly to the SP setting. Starting with initial value $x_1=\bar{x}_1$ and generated  sample path $\hat{\xi}_t$, $t=1,...,$
the policy $\bar{u}_t=\pi_t(x_t)$ is computed iteratively going forward in time   (compare with \eqref{conpol-2}). That is, at stage $t=1,...$, given value $\bar{x}_t$  of the state vector, the control value is computes as
\begin{equation}\label{conp-2a}
  \bar{u}_t \in  \argmin_{u_t\in \U_t(\bar{x}_t)}
\sum_{j=1}^N p_{tj} \left [c_{tj}(\bar{x}_t,u_t)+
\undv_{t+1}\big(A_{tj}\bar{x}_t +B_{tj}u_t+b_{tj}  \big)\right].
\end{equation}
Consequently  the next stage state value is computed
\begin{equation}\label{conp-2b}
\bar{x}_{t+1}= A_t (\hat{\xi}_t) \bar{x}_t +B_t (\hat{\xi}_t) \bar{u}_t +b_t (\hat{\xi}_t),
\end{equation}
and so on.
 Recall that in the SOC framework,  $\xi_1$ is also random.
 The computed values $\bar{x}_t$ of the state variables  can be used as trial points in the next iteration of the backward step of the algorithm.

\begin{remark}[$Q$-factor approach]
\label{rem-q}
{\rm
Consider  dynamic equations \eqref{soc-ind} with   affine state equations of the form \eqref{affine},   and
define\footnote{The $Q$-functions here should not be confused with the cost-to-go functions of the SP framework.}
\begin{equation}\label{qfa-1}
Q_t(x_t,u_t):=\bbe  \left[(c_{t}(x_t,u_t,\xi_t)+ V_{t+1} (A_{t} x_t+B_{t} u_t+ b_{t})\right].
\end{equation}
We have that
\begin{equation}
\label{qval}
 V_t(x_t)=
 \inf\limits_{u_t\in \U_t(x_t)}   Q_t(x_t,u_t).
 \end{equation}
 The dynamic equations  \eqref{soc-ind}  can be written in terms of $Q_t(x_t,u_t)$ as
\begin{equation}
\label{qfac-2}
Q_t(x_t,u_t)=\bbe \Big[c_{t}(x_t,u_t,\xi_t)+
 \inf\limits_{u_{t+1}\in \U_{t+1}(A_{t} x_t+B_{t} u_t+ b_{t})}   Q_{t+1}\big( A_{t} x_t+B_{t} u_t+ b_{t},   u_{t+1}\big)\Big].
\end{equation}
The
optimal policy $\bar{u}_t=\pi_t(x_t)$, defined by \eqref{conpol-2},  can be written in terms of the $Q$-functions as
\begin{equation}\label{qpoly}
\bar{u}_t\in \argmin_{u_t\in \U_t(x_t)} Q_t(x_t,u_t).
\end{equation}

Under Assumption \ref{ass-conv},   functions $Q_t(x_t,u_t)$ are convex  jointly in $x_t$ and $u_t$. Consequently   the cutting plane algorithm  can be applied  directly  to dynamic equations \eqref{qfac-2}. That is, in the backward step    current lower approximations $\underline{Q}_{t}(x_t,u_t)$ of the respective  functions
$Q_t(x_t,u_t)$ are updated by adding    cutting planes $\ell_t(x_t,u_t)$  computed at trial points $(\tilde{x}_t,\tilde{u}_t)$  generated in the forward step of of   the algorithm. The constructed
approximations $\underline{Q}_{t}(\cdot,\cdot)$ are used
in the forward step of the algorithm for generating  a point estimate of the expected value of the corresponding policy for a randomly generated sample path.
} $\hfill \square$
\end{remark}

\subsection{SDDP  algorithms for risk averse problems.}
 \label{sec-cut risk}

Cutting planes algorithms of
 SDDP type   can be extended to risk averse    problems  (cf., \cite{phil13},\cite{shaejor2011}). Let us start with the    SP framework.

\subsubsection  {Stochastic Programming framework.}
\label{sec-rasp}
We follow    the assumptions and notation  of Section \ref{sec-stprcut}, and assume that  equation \eqref{dynind-1a}  for the cost-to-go functions,   is  replaced by its risk averse counterpart  \eqref{drisk-1}. The backward step of the algorithm is performed similar to the risk neutral case. We only need to make adjustment to calculation of  subgradients of  the
cutting planes   approximations $\und_{t+1}(\cdot)$  of the
cost-to-go functions $\Q_{t+1} (\cdot)$.

In numerical algorithms, we deal with discrete distributions having a finite support. So we assume in this section that  $\O=\{\w_1,...,\w_N\}$ is  a finite space equipped with sigma algebra $\F$   of all subsets of $\O$,  and reference probability measure $\bbp$    with the respective {\em nonzero} probabilities $p_i=\bbp(\{\w_i\})$, $i=1,...,N$.
 For a random variable (a function) $Z:\O\to \bbr$ , we denote $Z_i:=Z(\w_i)$, $i=1,...,N$.

\begin{remark}
\label{rem-risk}
{\rm
Recall definition of law invariant risk measures
discussed in section \ref{sec-risk}. Law invariant coherent risk measures have naturally defined conditional counterparts used in construction of the corresponding nested  risk measures  (see equations \eqref{inf-risk} - \eqref{infr-1} below).
Suppose for the moment that the  probabilities $p_i$  are such that for $\I,\J\subset\{1,...,N\}$ the equality  $\sum_{i\in \I}p_i=\sum_{j\in \J}p_j$ holds only if $\I=\J$. Then two variables
$Z,Z':\O\to \bbr$ are distributionally equivalent iff they do coincide. Therefore in that case any risk measure defined on the space of variables $Z:\O\to \bbr$ is law invariant.
In the theory of law invariant coherent risk measures it is usually  assumed that the reference probability space $(\O,\F,\bbp)$ is  atomless; of course an atomless space cannot be finite.
In the case of finite space $\O$ it makes sense to consider law invariant risk measures   when all probabilities $p_i$ are equal to each other, i.e., $p_i=1/N$, $i=1,...,N$. In that case $Z,Z':\O\to \bbr$ are distributionally equivalent iff there  exists a permutation $\pi:\O\to\O$ such that $Z'_i=Z_{\pi(i)}$, $i=1,...,N$.
Anyway we will not restrict the following discussion to the case of  all probabilities $p_i$ being  equal to each other.
} $\hfill \square$
\end{remark}

 In the considered  setting
$\bbe_\bbp[Z]=\sum_{i=1}^N p_i  Z_i$, and
the dual representation \eqref{fmtheorem} of coherent risk measure
can be written as
\begin{equation}\label{finitedual}
\begin{array}{ll}
 \R(Z)=\sup\limits_{\zeta\in \cA} \sum_{i=1}^N p_i \zeta_i Z_i,
 \end{array}
\end{equation}
where    $\cA$ is a convex closed subset of $\bbr^N$  consisting of vectors
$\zeta\ge 0$ such that  $\sum_{i=1}^N p_i\zeta_i=1$.
Since it is  assumed  that $\R(Z)$ is real valued, the set $\cA$ is bounded and  hence is compact.

The subdifferential of $\R(\cdot)$ at $Z$ is given by
\begin{equation}\label{subdif}
\begin{array}{ll}
  \partial \R(Z)=\argmax\limits_{\zeta\in \cA}\sum_{i=1}^N p_i \zeta_i Z_i.
  \end{array}
\end{equation}
 Specific formulas for the subgradients, in various examples of coherent risk measures,   are listed in  \cite[Section 6.3.2]{SDR}.
These formulas can be applied for computation of cutting planes in the backward step of the algorithm.  Several such examples were discussed and implemented in an SDDP type  algorithm  in \cite{Sha2012a}.

As it was described in Section \ref{sec-stprcut}, in the risk neutral setting   the forward step of the algorithm has two functions; namely, construction of statistical upper bounds and generation of trial points. The constructed approximations $\und_{t+1}(x_t)$  of the cost-to-go functions define a feasible policy in the same way as in the risk neutral case by using formula \eqref{polappl}. Therefore the forward step can be used for generation  of trial points similar to  the risk neutral case. However, because of the nonlinearity of the risk measures,  it does not produce an unbiased estimate of the risk averse value of the constructed policy and cannot be applied in a straightforward way for computing  statistical upper bounds similar to the risk neutral case.  At the moment there are no known efficient numerical algorithms for computing  statistical upper bounds for risk averse problems  in the SP framework. Interestingly, as we will argue below, in the SOC setting it is possible to construct reasonably efficient statistical upper bounds for a certain class of risk measures.

\subsubsection{Stochastic Optimal Control framework.}
\label{sec-upper}
As before we suppose that  Assumptions \ref{assum-ind} and   \ref{assum-stage} hold, and thus value functions $V_t(x_t)$, associated with coherent risk measure $\R$, are determined by
the dynamic equations \eqref{drisk-3}. Furthermore, we suppose that Assumption \ref{ass-conv} is satisfied, and hence value functions $V_t(x_t)$ are convex.

We consider in this section the following class of coherent risk measures. This class is quite broad and includes  many   risk measures used in practical applications (the presentation below is based on    \cite{GSC2022}). We assume that the considered  risk measures can be represented as
\begin{equation}\label{riskmes}
 \R_t(Z):=\inf_{\theta\in \Theta}\bbe_{\bbp_t}[\Psi(Z,\theta)],\;\;Z\in \Z,
\end{equation}
where $\Theta$ is a nonempty convex closed subset of a finite dimensional vector space and $\Psi:\bbr\times\Theta\to \bbr$ is a real valued function. For the sake of simplicity, we assume that the function $\Psi(z,\theta)$ is the same at every stage, while the reference probability distribution $\bbp_t$  of $\xi_t$ can be different at different stages.
We make the following assumptions about function $\Psi(z,\theta)$.
(i)  For every $Z\in \Z$, the
expectation in the right hand side of \eqref{riskmes} is well defined and the infimum is finite valued.
(ii)  Function $\Psi(z,\theta)$ is convex in $(z,\theta)\in \bbr\times \Theta$.
  (iii) $\Psi(\cdot, \theta)$ is monotonically nondecreasing, i.e., if $z_1\le z_2$,  then $\Psi(z_1,\theta)\le \Psi(z_2,\theta)$ for all $\theta\in \Theta$.

 Under these assumptions, the  functional defined in \eqref{riskmes}  satisfies
 the axioms of convexity and monotonicity.
 Convex combination of expectation and   Average Value-at-Risk, defined in \eqref{convcomb},
    is of that  form  with $\Theta:=\bbr$,    $\Z:=L_1(\O,\F,\bbp)$ and
    \begin{equation}\label{psiform}
    \Psi(z,\theta):=(1-\lambda_t )z+ \lambda_t  \left(\theta+ (1-\alpha_t)^{-1}[z-\theta]_+\right).
    \end{equation}

For risk measures of the form \eqref{riskmes},  dynamic programming equations \eqref{drisk-3} can be written as
\begin{equation}
\label{drisk-fun}
V_t(x_t)=\inf\limits_{u_t\in \U_t(x_t),\,\theta\in \Theta}
\underbrace{\sum_{j=1}^N p_{tj} \Psi \big (c_{tj}(x_t,u_t)+
V_{t+1}(A_{tj}x_t +B_{tj}u_t+b_{tj} ),\theta\big)}_{\bbe_{\bbp_t}
[\Psi  (c_{t}(x_t,u_t,\xi_t)+
V_{t+1}(A_{t}x_t +B_{t}u_t+b_{t} ),\theta)]}.
\end{equation}
Note that  the minimization in the right hand side of \eqref{drisk-fun} is performed jointly in controls $u_t$ and parameter vector $\theta$.

The backward step of the algorithm can performed similar to the risk neutral case. In case   $\U_t(x_t)\equiv \U_t$ does not depend on $x_t$, by the chain rule of differentiation  formula \eqref{subd-3a} for the required  subgradient  is extended  to
\begin{equation}\label{subd-ab}
\sum_{j=1}^N  p_{tj}  \left [\Psi'(y_{tj},\bar{\theta}_t)\left(\nabla  c_{tj}(x_t,\bar{u}_t)+
A^\top_{tj}  \nabla \undv_{t+1}\big(A_{tj} x_t+B_{tj} \bar{u}_t+ b_{tj} \big)\right)\right],
\end{equation}
where $\bar{u}_t$ and $\bar{\theta}_t$ are minimizers in the   right hand side of \eqref{drisk-fun}, with $V_{t+1}(\cdot)$ replaced by its current approximation $\undv_{t+1}(\cdot)$,  and
$\Psi'(y_{tj},\bar{\theta}_t)$ is a subgradient of  $\Psi(\cdot,\bar{\theta}_t)$
at $y_{tj}:=c_{tj}(x_t,\bar{u}_t)+
\undv_{t+1}\big(A_{tj}x_t +B_{tj}\bar{u}_t+b_{tj}  \big)$.

The subgradient    $\Psi'(\cdot,\theta)$ can be easily computed in many interesting examples (e.g., \cite[Section 6.3.2]{SDR}).
For example,  consider convex combination of expectation and   Average Value-at-Risk measure, defined in \eqref{convcomb}.  The corresponding function $\Psi(z,\theta)$   of the form \eqref{psiform}, is piecewise linear in $z$. Then $\Psi'(z,\theta)=1-\lambda_t$ if $z<\theta$,
 $\Psi'(z,\theta)=1-\lambda_t+\lambda_t(1-\alpha_t)^{-1}$ if $z>\theta$; and if $z=\theta$,  then $\Psi'(z,\theta)$ can be any point in the interval
 $[1-\lambda_t, 1-\lambda_t+\lambda_t(1-\alpha_t)^{-1}]$.

The computed approximations $\undv_{t+1}(\cdot)$ of the value functions define a feasible policy in the same way as in the risk neutral case. Also by construction we have that $V_t(\cdot)\ge \undv_t(\cdot)$, $t=1,...,T$. Therefore the algorithm  produces a (deterministic) lower bound for the optimal value of the  problem.

Formula \eqref{drisk-fun} shows that the controls and value of the parameter $\theta$ can be computed simultaneously. This leads to the following statistical upper bound for the (risk averse) value of the policy defined by the computed approximation.   Let $\hat{\xi}_1,...,\hat{\xi}_T$ be a   sample path of the data process.
Suppose for the moment that $\R_t=\bbe$, i.e., consider the risk neutral case. Let $\pi$ be a policy determined by controls  $\bar{u}_t$ and   states $\bar{x}_t$, and  $\hat{c}_t:=c_t(\bar{x}_t,\bar{u}_t,\hat{\xi}_t)$, $t=1,...,T$, and $\hat{c}_{T+1}:= c_{T+1}(\bar{x}_{T+1})$, be cost values associated with the   sample path. Then   the corresponding   point estimate of the expectation \eqref{soc} of the  total cost for the considered policy $\pi$,   can be computed iteratively going backward in time with $\cv_{T+1}:=\hat{c}_{T+1}$ and
\begin{equation}\label{iterat}
  \cv_t:= \hat{c}_t+\cv_{t+1},\;\;t=T,...,1.
\end{equation}
 Value $\cv_1$ gives an unbiased estimate of the corresponding expected total cost.

This process can be adapted   for general risk averse case as follows.
Let $\hat{A}_t=A_t (\hat{\xi}_t)$, $\hat{B}_t=B_t (\hat{\xi}_t)$  and $\hat{b}_t=b_t (\hat{\xi}_t)$,  $t=1,...,T$,  be   realizations of the parameters corresponding to the generated sample path.
For the  approximations $\undv_t(\cdot)$ consider the corresponding policy with controls $\bar{u}_t$ and parameters $\bar{\theta}_t$ computed going forward,  starting with initial value $\bar{x}_1$ and using (compare with \eqref{conp-2a} - \eqref{conp-2b})
\begin{eqnarray}
\label{contpar-a}
 (\bar{u}_t,\bar{\theta}_t)&\in &
  \argmin\limits_{u_t\in \U_t(\bar{x}_t),\,\theta\in \Theta}
 \sum_{j=1}^N p_{tj} \Psi \big (c_{tj}(\bar{x}_t,u_t)+
\undv_{t+1}(A_{tj}\bar{x}_t +B_{tj}u_t+b_{tj} ),\theta\big), \\
\label{contpar-b}
 \bar{x}_{t+1}&=& \hat{A}_t  \ \bar{x}_t +\hat{B}_t  \bar{u}_t +\hat{b}_t.
\end{eqnarray}
Let $\hat{c}_t:=c_t(\bar{x}_t,\bar{u}_t,\hat{\xi}_t)$, $t=1,...,T$, and $\hat{c}_{T+1}:= c_{T+1}(\bar{x}_{T+1})$, be cost values associated with the generated sample path.
Consider the following values defined iteratively going backward in time: $\cv_{T+1}:=\hat{c}_{T+1}$ and
\begin{equation}\label{iterat2}
  \cv_t:=\Psi(\hat{c}_t+\cv_{t+1},\bar{\theta}_t),\;\;t=T,...,1.
\end{equation}

It is possible to show that $\bbe[\cv_1]$  is greater than or equal to   value of the policy defined by the considered
approximate value functions, and hence  is an upper bound on the optimal value of the risk averse
problem (cf., \cite{GSC2022}).
Expectation  $\bbe[\cv_1]$ can be estimated by randomly generating  sample paths (scenarios) and averaging the computed  realizations of random variable  $\cv_1$. This is similar to construction of statistical upper bound in the risk neutral setting. The inequality:
``$\bbe[\cv_1]\ge$ value of the constructed policy",
is based on Jensen’s inequality applied to convex function $\Psi(\cdot,\theta)$ and, unlike the risk neutral case,  can be strict. Therefore  the above procedure  introduces  an additional error in the estimated optimality gap.   This additional error becomes larger when   function $\Psi(\cdot,\theta)$  is ``more nonlinear". Numerical experiments indicate that this procedure gives a reasonable upper bound, for instance  for risk measures  of the form \eqref{convcomb}  (cf., \cite{GSC2022}).

\paragraph{$Q$-factor approach. }
Consider the dynamic equations \eqref{drisk-fun} and define (see Remark \ref{rem-q})
\begin{equation}\label{qfact-1}
Q_t(x_t,u_t,\theta_t):=\bbe_{\bbp_t} \left[ \Psi\big (c_{t}(x_t,u_t,\xi_t)+ V_{t+1} (A_{t} x_t+B_{t} u_t+ b_{t}),\theta_t  \big )\right].
\end{equation}
We have that
\[
 V_t(x_t)=
 \inf\limits_{u_t\in \U_t,\,\theta_t\in \Theta}   Q_t(x_t,u_t,\theta_t),
 \]
 and hence   dynamic equations  \eqref{drisk-fun}  can be written in terms of $Q_t(x_t,u_t,\theta_t)$ as
\begin{equation*}
Q_t(x_t,u_t,\theta_t)=\bbe_{\bbp_t} \left[ \Psi  \Big(c_{t}(x_t,u_t,\xi_t)+
 \inf\limits_{u_{t+1}\in \U_{t+1},\,\theta_{t+1}\in \Theta}   Q_{t+1}\big( A_{t} x_t+B_{t} u_t+ b_{t},   u_{t+1},\theta_{t+1} \big),\theta_t\Big)\right].
\end{equation*}

The cutting planes   algorithm can be applied directly to  functions $Q_t(x_t,u_t,\theta_t)$ rather than to the value functions $V_t(x_t)$. In the backward step of the algorithm,  subgradients with respect to $x_t,u_t$ and $\theta_t$,   of the current approximations of  $Q_t(x_t,u_t,\theta_t)$,  should be computed at the respective trial points,
and the obtained cutting plane $\ell_t(x_t,u_t,\theta_t)$ is added to the current family of cutting planes of $Q_t(x_t,u_t,\theta_t)$.
An advantage of that approach could be  that the calculation of the respective subgradients does not require solving {\em nonlinear optimization} programs  even if the function $\Psi$ is not polyhedral.

\subsubsection{Infinite horizon setting}
\label{sec-inf}
 An infinite horizon (stationary)   counterpart of the SOC problem  \eqref{soc} - \eqref{soc-b}  is
 \begin{equation}\label{inf-soc}
\min\limits_{\pi\in \Pi} \bbe^\pi\Big [ \sum_{t=1}^{\infty}
\gamma^{t-1}c(x_t,u_t,\xi_t)
\Big],
 \end{equation}
 where $\gamma\in (0,1)$   is the so-called  discount factor.
 The optimization (minimization) in \eqref{inf-soc} is over the set of feasible  policies $\Pi$ satisfying\footnote{When the number of scenarios is finite, the feasibility constrains should be satisfied for all possible realizations of the random process $\xi_t$.} (w.p.1):
  \begin{equation}\label{infsoc-2}
u_t\in \U,\;x_{t+1}=F(x_t,u_t,\xi_t),\;t\ge 1.
 \end{equation}
 It is assumed that $\xi_t\sim P$, $t=1,...,$  is an i.i.d. (independent identically distributed)  sequence  of random vectors, with common distribution $P$. The cost $c(x,u,\xi)$ and the mapping  $F(x,u,\xi)$  are assumed to be  the same at every stage, and the (nonempty)  feasibility set $\U$ of controls  does not depend on the state and stage.

 Bellman equation for the value function, associated with problem \eqref{inf-soc}, can be written as
\begin{equation}\label{bell-1}
 V(x)=\inf_{u\in \U}\bbe_P[c(x,u,\xi) +\gamma V(F(x,u,\xi))].
\end{equation}
 Assuming that
 $c(x,u,\xi)$ is bounded,  equation \eqref{bell-1} has a unique solution
$\bar{V}(\cdot)$. The optimal policy for problem \eqref{infsoc-2}  is given by $u_t=\pi(x_t)$, $t=1,...,$   with
\begin{equation}\label{bell-2}
\pi(x)\in \argmin_{u\in \U}\bbe_P\big[c(x,u,\xi) +\gamma \bar{V}(F(x,u,\xi))\big].
\end{equation}

Similar to Assumption \ref{ass-conv} let us make following assumption ensuring convexity of the solution $\bar{V}(\cdot)$ of equation \eqref{bell-1}.

\begin{assumption}
The function $c(x,u,\xi)$ is convex in $(x,u)$, the mapping $F(x,u,\xi)=A(\xi)x+B(\xi)u+b(\xi)$ is affine, and the set $\U$ is convex closed.
\end{assumption}

In order to solve equation \eqref{bell-1} numerically we need to discretize the possibly continuous distribution $P$. Suppose that the SAA method is applied, i.e., a sample  $\xi^1,...,\xi^N$ of size $N$ from the distribution $P$ is generated (say by Monte Carlo sampling techniques), and the distribution $P$ in Bellman equation \eqref{bell-1} is replaced by its empirical estimate\footnote{By $\delta_\xi$ we denote measure assigning mass one at point $\xi$, the so-called Dirac measure.} $P_N=N^{-1}\sum_{j=1}^N \delta_{\xi^j}$, assigning probability $1/N$ to each generated point. This raises the question of the involved sample complexity, i.e., how large should be the
sample size $N$ in order for the SAA problem to give an accurate approximation of the original
problem. In some applications the discount factor
 is very close to one. It is well known
that as the discount factor $\gamma$ approaches one, it becomes more difficult to solve the problem. Note that the optimal value of problem \eqref{inf-soc}
increases at the rate of $O((1-\gamma)^{-1})$ as $\gamma$
 approaches one. It is shown  in \cite{shacheng2021} that the required  sample size $N$ is  of order $O((1-\gamma)^{-1}\e^{-2})$ as a function
of the discount factor $\gamma$
 and the error level $\e>0$. That is, the error grows more or less  proportionally to the optimal value as $\gamma$ approaches one. This indicates that the sample size required to achieve a specified {\em relative} error  is not   sensitive to the discount factor being close to one.

Suppose now that random vector $\xi$ has a finite number $\xi^1,...,\xi^N$ of possible values with respective probabilities $p_1,...,p_N$. For example, in the SAA approach this is  a random sample generated by Monte Carlo sampling techniques, in which case $p_i=1/N$, $i=1,...,N$.   In that setting  Bellman equation \eqref{bell-1} can be written as
\begin{equation}\label{bell-dis}
 V(x)=\inf_{u\in \U}\, \sum_{j=1}^N p_j[c_j(x,u) +\gamma V(A_jx+B_ju+b_j)],
\end{equation}
where $c_j(x,u):=c(x,u,\xi^j)$, $A_j:=A(\xi^j)$, $B_j:=B(\xi^j)$ and $b_j:=b(\xi^j)$, $j=1,...,N$.

A cutting planes algorithm of SDDP type  can be applied to equation \eqref{bell-dis}. Given a current approximation $\undv (\cdot)$ of the value function by cutting planes and a trial point $\tilde{x}$, new cutting plane $\ell(x)=\tilde{v}+\sg^\top (x-\tilde{x})$ is added to the set of cutting planes of the value function,  where
\begin{eqnarray}
\label{bell-d3}
\tilde{u}\in&\argmin\limits_{u\in \U}\sum_{j=1}^N p_j[c_j(\tilde{x},u) +\gamma \undv(A_j\tilde{x}+B_ju+b_j)],\\
\tilde{v} =&\sum_{j=1}^N p_j[c_j(\tilde{x},\tilde{u}) +\gamma \undv(A_j\tilde{x}+B_j\tilde{u}+b_j)],\\
\sg=&\sum_{j=1}^N p_j[\nabla c_j(\tilde{x},\tilde{u}) +\gamma A_j^\top \nabla \undv(A_j\tilde{x}+B_j\tilde{u}+b_j)].
\label{bell-d2}
\end{eqnarray}
 If $c(x,u,\xi)$ is linear (polyhedral)  in $u$ and the set $\U$ is polyhedral given by linear inequalities,  the minimization problem in the right hand side of \eqref{bell-d3} can be formulated as a linear programming problem.  The subgradient of $\undv(\cdot)$  at a point $x$  can be computed in the same way as in
\eqref{subd-3a} by using subgradient of its supporting  plane  at $x$.

Similar to the case of finite horizon,  by construction we have    that $\bar{V}(\cdot) \ge \undv(\cdot)$, and hence $\undv(x_1)$ gives a (deterministic) lower bound for the optimal value of problem \eqref{inf-soc} with initial value $x_1$  of the state vector.
Important questions, related to the above  algorithm,  is how to generate trial points and  how to construct an upper bound for the  optimal value of problem \eqref{inf-soc}.
Current implementation, which was tested in numerical experiments, is to apply the forward step of the cutting planes algorithm to a truncated version of problem \eqref{inf-soc}. Note that
\begin{equation}\label{geom}
 \left|\sum_{t=T+1}^{\infty}
\gamma^{t-1}c(x_t,u_t,\xi_t)\right|\le
\sum_{t=T+1}^{\infty}
\gamma^{t-1}\left| c(x_t,u_t,\xi_t)\right|\le
\frac{\kappa \gamma^T}{1-\gamma},
\end{equation}
where $\kappa$ is an upper bound for $|c(x,u,\xi)|$.  Therefore the suggestion is to choose the horizon $T$ such that $\kappa \gamma^T(1-\gamma)^{-1}\le \e$, where $\e>0$  is     a prescribed accuracy.  Then  the forward steps are applied  to the finite horizon version of the problem using current approximation $\undv(\cdot)$ of the value function. One run of the forward step, applied to the truncated version, will generate a point estimate (up to accuracy $\e$) of   value of the  constructed policy,  and $T-1$ trial points. It is possible to choose, say at random, one of these trial points    for construction of the cutting plane for the value function.

Such choice  of the  trial points could be not optimal. The question of ``best way" for construction of trial points in the above algorithm  was not carefully investigated. Another problem with this procedure is that when the discount factor is close to one, the required   horizon $T$ could be too large for a numerical implementation; recall that in order to construct the statistical upper bound the forward step  should be run a reasonable  number of times.
In such cases,  if the problem can be formulated as a linear program, deterministic upper bounds of the dual problem can be more efficient
(cf., \cite{shacheng2021a}).

\paragraph{Risk averse case.} The infinite horizon  risk neutral formulation \eqref{inf-soc} can be extended to the risk averse setting by
replacing the expectation operator  in Bellman equation \eqref{bell-1} with  a    coherent risk measure $\R$. That is, risk averse counterpart of Bellman equation \eqref{bell-1} is
\begin{equation}\label{bellrisk}
 V(x)=\inf_{u\in \U}\R[c(x,u,\xi) +\gamma V(F(x,u,\xi)].
\end{equation}
This equation corresponds to the nested counterpart of the risk neutral problem \eqref{inf-soc}:
 \begin{equation}\label{inf-risk}
 \begin{array}{ll}
\min\limits_{\pi\in \Pi}\ \lim\limits_{T\to\infty} \cR_T\left ( \sum_{t=1}^{T}
\gamma^{t-1}c(x_t,u_t,\xi_t)
\right),
\end{array}
\end{equation}
where  $ \cR_T$ is the corresponding $T$-stage nested risk measure (cf., \cite{RuSh:2006b}).  That is, as it was discussed in section \ref{sec-soc},
for a considered policy, state $x_t$ and control $u_t$ at stage $t$ are functions of $\xi_{[t-1]}$. Consider the corresponding cost $Z_t:=\gamma^{t-1}c(x_t,u_t,\xi_t)$ which is a function of $\xi_{[t]}$. Then
 \begin{equation}\label{infr-1}
 \begin{array}{ll}
\cR_T\left ( \sum_{t=1}^{T}
Z_t
\right)=\R\Big(Z_1+\R_{|\xi_{[1]}}\big(Z_2+...
+\R_{|\xi_{T-1}}(Z_{T})\big)\Big),
\end{array}
\end{equation}
where $\R_{|\xi_{[t]}}$ is the conditional counterpart of $\R$.

The backward step of the  cutting planes algorithm in the risk averse setting is basically the same as the one discussed above for the risk neutral case. In order to construct a statistical upper bound we can proceed similar to the construction of section \ref{sec-upper}. That is, suppose that risk measure $\R$ is of the form \eqref{riskmes} with function $\Psi(z,\theta)$ satisfying conditions (i) - (iii) specified in section
\ref{sec-upper}. Suppose also  that $\xi$ has a finite number $\xi^1,...,\xi^N$  of possible values with
reference  probabilities $p_1,...,p_N$.
Then  Bellman equation \eqref{bellrisk}  can be written as
\begin{equation}\label{bellris-2}
V(x)=\inf_{u\in \U,\,\theta\in \Theta}\, \sum_{j=1}^N p_j\Psi\big(c_j(x,u) +\gamma V(A_jx+B_ju+b_j),\theta\big).
\end{equation}
A statistical upper bound can be constructed in a way similar to section \ref{sec-upper} using a truncated version of problem \eqref{inf-soc}.

\paragraph{Periodical setting}
 In various applications the data process $\xi_t$ has a periodical behavior. One such example is hydropower generation planning  discussed in  \cite{Sha2012a}. The uncertain process in that example is water inflows recorded  on monthly basis. There is available
 historical data    of 79 years of   records of
   the natural monthly energy inflow. This allows to model the uncertain process as a periodical random process with period $\sm=12$. The planning in that example  was for 5 years with another 5 years to eliminate  the   end of  horizon effect. This resulted in $T=120$ stage stochastic optimization problem. Periodical approach, which we discuss below, allows drastic reduction in the number of stages while giving almost the same policy value  for the first stage.

 In \cite{shading} such periodical approach was introduced from the stochastic programming point of view.
 We discuss this  in the framework of   SOC. Consider the following infinite horizon SOC problem with discount factor $\gamma\in (0,1)$:
 \begin{eqnarray}
 \label{eqperiod-1}
 \min\limits_{u_t\in \U_t} & \bbe^\pi\left [ \sum\limits_{t=1}^{\infty}
\gamma^{t-1}c_t(x_t,u_t,\xi_t) \right] \\
 \label{eqperiod-2}
{\rm s.t.} & x_{t+1}=F_t(x_t,u_t,\xi_t),\; t\ge 1.
 \end{eqnarray}
 We assume   the following  periodic structure, with period $\sm\ge 1$:
 the random data  process $\xi_t$ is stagewise independent   with  $\xi_t$ and $\xi_{t+\sm}$ having  the same probability distribution, and with
 $c_t(\cdot,\cdot,\cdot)=c_{t+\sm}(\cdot,\cdot,\cdot)$,
 $F_t(\cdot,\cdot,\cdot)=F_{t+m}(\cdot,\cdot,\cdot)$  and $\U_t=\U_{t+\sm}$, for $t\ge 1$.

Bellman equations for this problem take the form  (compare with \eqref{soc-ind})

\begin{eqnarray}
\label{soc-p}
 V_t(x_t)&=&\inf\limits_{u_t\in \U_t}
\bbe \left [c_t(x_t,u_t,\xi_t)+ \gamma
V_{t+1}\big(F_t(x_t,u_t,\xi_t) \big)\right],\;t=1,...,\sm-1,\\
\label{soc-per}
V_\sm(x_\sm)&=&\inf\limits_{u_\sm\in \U_\sm}
\bbe \left [c_\sm(x_\sm,u_\sm,\xi_\sm)+ \gamma
V_{1}\big(F_\sm(x_\sm,u_\sm,\xi_\sm) \big)\right],\;t=\sm.
   \end{eqnarray}
For period $\sm=1$,   periodical problem \eqref{eqperiod-1} -\eqref{eqperiod-2}
coincides with   the stationary problem \eqref{inf-soc}
 of section \ref{sec-inf}, and equation \eqref{soc-per}  becomes Bellman equation \eqref{bell-1}.
 The optimal
policy for problem \eqref{eqperiod-1} - \eqref{eqperiod-2}  is given   by  $u_t=\pi_t(x_t)$  with
\begin{eqnarray}
\label{perpol-1}
  \pi_t (x_t)&\in&\inf\limits_{u_t\in \U_t}
\bbe \left [c_t(x_t,u_t,\xi_t)+\gamma
\bar{V}_{t+1}\big(F_t(x_t,u_t,\xi_t) \big)\right],\;t=1,...,\sm-1,\\
\label{perpol-2}
 \pi_\sm(x_\sm)&\in&\inf\limits_{u_\sm\in \U_\sm}
\bbe \left [c_\sm(x_\sm,u_\sm,\xi_\sm)+ \gamma
\bar{V}_{1}\big(F_\sm(x_\sm,u_\sm,\xi_\sm) \big)\right],\;t=\sm,
   \end{eqnarray}
and $\pi_{t+\sm}(\cdot)=\pi_t(\cdot)$ for $t\ge 1$, where  $(\bar{V}_1,...,\bar{V}_\sm)$ is the solution  of equations \eqref{soc-p} - \eqref{soc-per}.

Such periodical formulation can be extended to the (nested) risk averse setting and
cutting planes type algorithm can be applied in a way similar to what was discussed above  (cf., \cite{shading}).

  \setcounter{equation}{0}
\section{Computational Complexity of Cutting Plane Methods}

As mentioned earlier, SDDP reduces to Kelly's cutting plane method for two stage
problems and it can behave poorly as the number of decision variables increases.
On the other hand, numerical experiments indicate  that
this type of algorithm can scale well with respect to  the number of stages $T$.
In this  section  we summarize some recent theoretical studies in~\cite{Lan22,JuLan23} on the rate of
convergence of SDDP. We first establish the number of iterations required by the SDDP method, and then present a variant of SDDP that can improve
this complexity bound in terms of the dependence on the number of stages. Related developments, with a focus more on mixed-integer nonlinear optimization,
can also be found in \cite{ZhangSun22-1}.

\subsection{Computational complexity of SDDP}
Consider the multi-stage stochastic programming  problem  \eqnok{stpr-1}-\eqnok{stpr-1a}.
%
Similar to the the previous section,
we assume that $\xi_t$'s are stagewise independent (i.e., Assumption~\ref{assum-stage}) and
 finitely supported (i.e., Assumption~\ref{ass-finite}).
For simplicity we further assume that the probabilities for $\xi_{ti}$, $i = 1, \ldots, N$, are equal, i.e.,
$p_{ti}=1/N$, $t=2,...,T$,  $i=1,...,N$.
In that  case, the dynamic equations in \eqnok{dynind-1}-\eqnok{dynind-1a}
reduce to
\begin{equation}\label{define_value_function_sddp_random}
Q_t(x_{t-1},\xi_{ti})=  \inf_{x_t\in \X_t} \left\{ F_{ti}(x_t):= f_t(x_t,\xi_{ti})+\Q_{t+1}(x_{t}): B_{ti} x_{t-1} +A_{ti} x_t=b_{ti}\right\},
\end{equation}
where
\begin{equation}\label{define_value_function_sddp}
  \Q_{t+1}(x_{t})= \tfrac{1}{N} \tsum_{i=1}^N Q_{t+1}(x_t,\xi_{t+1,i}).
\end{equation}


The basic scheme of the SDDP method has been discussed in Section \ref{sec-stprcut}.
Here we add a little more details to facilitate its complexity analysis.
Each iteration of the SDDP method to solve problem \eqnok{stpr-1}-\eqnok{stpr-1a} contains
a forward  step and a backward  step. In the forward  step of the SDDP method,
we randomly pick up an index $i_t$ out of $\{1, \ldots, N\}$ and solve problem
\[
\begin{array} {lll}
x_{t}^k &\in & \argmin\left\{ \underline F_{t{i_t}}^{k-1}(x_t) := f_t(x_t, \xi_{t{i_t}}) + \und_{t+1}^{k-1} (x_t)\right\} \\
&& \ \ \ {\mbox{s.t.} \ \ \  B_{t i_t} x_{t-1}^k +A_{ti_t} x_t=b_{ti_t}},\; x_t \in \X_t,
\end{array}
\]
to update $x_t^k$. Equivalently, one can view $x_t^k$ as being randomly chosen from $\tilde x_{t i}^k$, $i = 1, \ldots, N$,
defined as
\beq \label{def_x_k_t_sddp}
\begin{array} {lll}
\tilde x_{ti}^k &\in& \argmin\left\{ \underline F_{ti}^{k-1}(x_t) := f_t(x_t, \xi_{ti}) +  \und_{t+1}^{k-1} (x_t)\right\}\\
&& \ \ \ {\mbox{s.t.} \ \ \  B_{t i} x_{t-1}^k +A_{ti} x_t=b_{ti}}, \;x_t \in \X_t.
\end{array}
\eeq
Note that we do not need to compute $\tilde x_{t i}^k$ for $i \neq i_t$, even though they will be used in the analysis of
the SDDP method. 
In next section, we will present a deterministic dual dynamic programming
method which chooses the feasible solution $x_{t}^k$ in the forward  step in a more
aggressive manner.

In the backward  step of SDDP, starting from the last stage, we update the cutting plane models
$\und_t^k$ according to
\[
\und_{t}^k(x_{t-1}) = \max \left \{\und_{t}^{k-1}(x),
\tfrac{1}{N} \tsum_{i=1}^{N} \left[ \cQ_{ti}^k(x_{t-1}^k)  + \langle (\cQ_{ti}^k)'(x_{t-1}^k),
x_{t-1} - x_{t-1}^k \rangle  \right]  \right\},
\]
for $t = T, T-1, \ldots, 2$.
Here $\cQ_{ti}^k(x_{t-1}^k)$ and $(\cQ_{ti}^k)'(x_{t-1}^k)$, respectively, denote the optimal value
and a subgradient at the point $x_{t-1}^k$ for the following approximate value function
\[
\begin{array} {lll}
\cQ_{ti}^k(x_{t-1}^k) &=& \min \left\{ \underline F_{ti}^{k}(x_t)  := f_t(x_t, \xi_{ti}) +  \und_{t+1}^{k} (x_t) \right\}\\
&& \ \ {\mbox{s.t.} \ \  B_{t i} x_{t-1}^k +A_{ti} x_t=b_{ti}}, x_t \in \X_t.
\end{array}
\]

In order to assess the progress made by each SDDP iteration, we need to
introduce a few notions.
We
say that a search point $x_{t}^k$ is $\epsilon_t$-{\sl saturated}
at iteration $k$ if
$
\Q_{t+1}(x_t^k) - \und_{t+1}^k(x_t^k) \le \epsilon_t.
$
Intuitively, a point $x_t^k$ is saturated if we have a good approximation
for the value function  $\cQ_{t+1}$ at $x_t^k$.
Moreover,
we say an $\epsilon_t$-saturated search point $x_t^k$ at stage $t$
is $\delta_t$-distinguishable if
$
\| x_t^k - x_t^j\| > \delta_t
$
for all other $\epsilon_t$-saturated search points $x_t^j$, $j \le k-1$, that have been generated
for stage $t$ so far by the algorithm.
Let us denote
$S_t^{k-1}$ the set of saturated points at stage $t$.
An $\epsilon_t$-saturated search point $x_t^k$
is $\delta_t$-distinguishable if
$
\dist(x_t^k, S_t^{k-1}) > \delta_t.
$
Recall that  $\dist(x, S)$  denotes the distance from point  $x$
to the set $S$.
Hence, in this case, this search point $x_t^k$ is far away enough from other saturated points.
The rate of convergence of SDDP depends on
how fast it generates $\epsilon_t$-saturated and $\delta_t$-distinguishable
search points. For simplicity we assume in this paper that
$\epsilon_t = \delta_t = \epsilon,\; t = 1, \ldots, T$,  for a given tolerance $\epsilon > 0$.

Under certain regularity assumptions, it  can be shown  that the function $F_{ti}$ and lower approximation functions $\underline F_{ti}^k$
are Lipschitz continuous with constant $M$.
Utilizing these continuity assumptions,
it  can be  shown  that
the probability for SDDP to
find a new $\epsilon$-distinguishable and $\epsilon$-saturated  search point at
the $k$-iteration can be bounded from below by
\beq \label{bnd_iter_prob_0}
 \tfrac{1}{\bar N} (1- \prob\{\tilde g_{t}^k \le \epsilon, t=1, \ldots, T-1\}),
\eeq
where
\beq \label{def_max_N_t}
\bar N :=  N^{T-2} \ \ \mbox{and} \ \ \
\tilde g_t^k :=  \tfrac{1}{N}\tsum_{i=1}^{N} \dist(\tilde x_{ti}^k, S_t^{k-1}).
\eeq
On the other hand, if for some iteration $k$, we have $\tilde g_t^k \le \epsilon$ for all $t = 1, \ldots, T-1$,
then
\beq \label{gap_bnd_ass1_sddp0}
F_{11}(x_1^k) - F^* \le F_{11}(x_1^k) - \underline F^{k-1}_{11}(x_1^k) \le 2 M (T-1) \, \epsilon.
\eeq
Note that the upper bound $2 M (T-1)$ can be further reduced so that it does not depend on $T$
if a discounting factor is incorporated in each stage of the problem.

Now let $K$ denote the number of iterations performed by the SDDP method before it
finds a forward path $(x_1^k, \ldots, x_T^k)$  such that \eqnok{gap_bnd_ass1_sddp0} holds.
Using the previous observation in \eqnok{bnd_iter_prob_0},
it  can be  shown  (see \cite{Lan22}) that $\bbe[K] \le  \bar K_\epsilon \bar N  + 2$, where $\bar N$ is defined in \eqnok{def_max_N_t} and
\beq \label{def_N_s_sddp}
\bar K_\epsilon := (T-1)\left( \tfrac{D}{\epsilon} +1 \right)^{n}.
\eeq
Here $n_t \le n$ and $D_t \le D$ for $t=1, \ldots, T$ denote
the dimension and diameter of the feasible sets $\X_t$, respectively.
In addition, for any $\alpha \ge 1$, we have
\[
\prob\{K \ge \alpha \bar K_\epsilon  \bar N + 1\} \le  \exp\left(-\tfrac{ (\alpha-1)^2 \bar K_\epsilon^2}{2 \alpha \bar N} \right).
\]

We now add a few remarks about the above complexity results for SDDP. Firstly, since SDDP is a randomized
algorithm, we provide bounds on the expected number of iterations required to find an approximate
solution of problem~\eqnok{stpr-1}-\eqnok{stpr-1a}. We also show that the probability of having
large deviations
from these expected bounds for SDDP decays exponentially fast.
Secondly, the complexity bound for
the SDDP method depends on $\bar N$, which increases exponentially with respect to $T$.
We will show in the next  section how to reduce the dependence on $T$ by presenting a variant of
the SDDP method.

\subsection{Explorative Dual Dynamic Programming}
The SDDP method in the previous section chooses the feasible solution $x_t^k$ in
a randomized manner.
In this section, we will introduce a deterministic method, called Explorative Dual Dynamic Programming (EDDP)
first presented in \cite{Lan22},
which chooses the feasible solution in the forward  step in an
aggressive manner.  As we will see, the latter approach will exhibit better
iteration complexity while the former one is easier to implement.

The EDDP method consists of the forward  step and backward  step.
In the forward  step of EDDP, for each stage $t$, we solve $N$ subproblems defined
in \eqnok{def_x_k_t_sddp}
to compute the search points $\tilde x_{ti}^k$, $i = 1, \ldots, N$.
For each  $\tilde x_{ti}^k$, we further compute the quantity $\dist(\tilde x_{ti}^k, S_t^{k-1})$,
the distance between $\tilde x_{ti}^k$ and the set $S_t^{k-1}$ of currently saturated search points in stage $t$.
Then we will choose
from $\tilde x_{ti}^k$, $i = 1, \ldots, N$, the one with the largest value of $\dist(\tilde x_{ti}^k, S_t^{k-1})$
as $x_t^k$, i.e.,
$\dist(x_t^k, S_t^{k-1}) = \max_{i=1, \ldots, N} \dist(\tilde x_{ti}^k, S_t^{k-1})$.
We can break the ties arbitrarily.
In view of the above discussion, the EDDP method always chooses the most ``distinguishable" forward path to encourage exploration
in an aggressive manner. This also explains the origin
of the name EDDP.
The backward  step of EDDP is similar to SDDP.

Under the same regularity assumptions as in SDDP,
it can be shown that each EDDP iteration will either
 generate at least one $\epsilon$-distinguishable and
$\epsilon$-saturated search point at some stage $t = 1, \ldots, T$,
or
find a feasible solution $x_1^k$ of problem \eqnok{stpr-1}-\eqnok{stpr-1a} such that
\begin{align}
F_{11}(x_1^k) - F^* &\le 2 M (T-1) \, \epsilon. \label{eq:bndOptGap_sddp}
\end{align}
As a result, the total number of iterations performed by EDDP before finding a feasible policy of problem \eqnok{stpr-1}-\eqnok{stpr-1a}
satisfying \eqnok{eq:bndOptGap_sddp} within at most $\bar K_\epsilon+1$ iterations where $\bar K_\epsilon$ is defined
in \eqnok{def_N_s_sddp} (see \cite{Lan22}). This complexity bound does not depend on the number of scenarios $\bar N$, and hence significantly
outperforms the one for the original SDDP method.
Moreover, we can show that the relation $\dist(x_1^k, S_1^{k-1}) \le \delta$ implies the relation in \eqnok{eq:bndOptGap_sddp}
for a properly chosen value of $\delta$,
and hence it can be used as a termination criterion for the EDDP method.
Note that it is possible to further reduce the dependence
 of the computational complexity of EDDP on $T$ so that the term $\bar K_\epsilon$ in \eqnok{def_N_s_sddp} does not
depend on $T$. However, this would require us to further modify EDDP by incorporating early termination
of the forward  step. More specifically,  we need to
terminate the forward  step at some stage $\bar t \le T$ and start the backward  step from
$\bar t$, whenever $\dist(x_{\bar t}^k, S_{\bar t}^{k-1})$ falls within a certain threshold value (see \cite{ZhangSun22-1,JuLan23}).

It should be noted that although the complexity of SDDP is worse than that for EDDP,
its performance in earlier  steps of the algorithm should be similar to that of EDDP.
Intuitively, for earlier iterations, the tolerance parameter $\delta_t$ used to define $\delta_t$-distinguishable
points are large.
As long as $\delta_t$ are large enough so that the solutions $\tilde x^k_{t i}$ are contained
within a ball with diameter roughly in the order of $\delta_t$, one can choose any
point randomly from $\tilde x^k_{t i}$ as $x^k_t$. In this case, SDDP will perform similarly to EDDP. This may explain why
SDDP exhibits good practical performance for low accuracy region.
For high accuracy region, the new EDDP algorithm
seems to be a much better choice in terms of its theoretical complexity.

\subsection{Complexity of EDDP and SDDP over an infinite horizon}

It is possible to generalize EDDP and SDDP to solve stochastic programs over an infinite horizon
(see Section~\ref{sec-inf}),
when the number of stages $T = \infty$, the cost function $f_t = f$,
the feasible set $\X_t = \X$ for all $t$, and the random variables $\xi_t$
are iid (independent identically distributed) and can be viewed as realizations
 of random vector  $\xi$.
A common line is to approximate the infinite horizon with a finite horizon. For example,
using the relation \eqnok{geom} in the SOC problem,
a horizon length of $T
=
O((\epsilon(1-\lambda))^{-1})$ suffices for a target accuracy $\epsilon$.
A direct extension of EDDP method show that
the computational complexity depends only quadratically on the time horizon $T$.
Exploiting the fact that  the infinite-horizon problem is stationary,
it has been shown that a simple modification that combines the forward and backward  step can
 improve the dependence on $T$ from quadratic to linear. Alternative selection strategies, including the use of upper bounds on the cost-to-go function and random sampling
as in SDDP can also be incorporated (see \cite{JuLan23} for more details).

 \setcounter{equation}{0}
\section{Dynamic Stochastic Approximation Algorithms}
\label{sec-sa}

\subsection{Extension of Stochastic Approximation}

Stochastic approximation (SA) has attracted much attention recently
for solving static stochastic optimization problems
given in the form of
\begin{equation} \label{sp}
\min_{x \in \X} \left\{f(x) := \bbe_\xi[F(x, \xi)] \right\},
\end{equation}
where $\X$ is a closed convex set, $\xi$ denotes the random vector and $F(\cdot, \xi)$
is a lower semicontinuous convex function.
Observe that we can cast two-stage stochastic optimization in the form
of \eqnok{sp} by viewing $F$ as the value function of the second-stage problem (see \cite{NJLS09-1,lns11,Lan10-3}).

The basic SA algorithm, initially proposed by Robbins and Monro~\cite{RobMon51-1}, mimics the simple projected gradient descent method by replacing exact
gradient with its unbiased estimator.
Important improvements for the SA methods have been made by Nemirovski and Yudin~\cite{nemyud:83}
and later by Polayk and Juditsky \cite{pol90,pol92}.  During the past few years, significant progress has been made
in SA methods including the incorporation of momentum and variance reduction, and  generalization to nonconvex setting
\cite{LanBook2020}.
It has been shown in \cite{NJLS09-1,lns11} that SA type methods can significantly outperform
the SAA approach for solving static (or two-stage) stochastic programming problems.
However, it remains unclear whether these SA methods can be generalized for multi-stage stochastic optimization
problems with $T \ge 3$.

In this section, we discuss a recently developed dynamic stochastic approximation (DSA) method for multi-stage stochastic optimization \cite{LanZhou17-1}.
The basic idea of the DSA method is to apply an inexact primal-dual SA method
for solving the $t$-th stage optimization problem to
compute an approximate stochastic subgradient for its associated value functions.
For simplicity, we focus on the basic scheme of the DSA algorithm and its main convergence properties
for solving three-stage stochastic optimization problems.
We also briefly discuss how to generalize it for solving more general form of multi-stage stochastic optimization with $T > 3$.

The main problem of interest in this section is the following three-stage SP problem:
\begin{equation}\label{3-stage prob}
\begin{aligned}
\min \, & f_1(x_1,\xi_1)+ &\bbe_{|\xi_1} [\min \, & \ f_2(x_2,\xi_2)& +\bbe_{|\xi_{[2]}}[\min \, & \ f_3(x_3, \xi_3)]]\\
  \text{s.t.} &\ A_1 x_1 = b_1,   & \text{s.t.} &\ A_2 x_2  +B_2 x_1 =  b_2, &  \text{s.t.}&\ A_3 x_3 + B_3x_2 =  b_3, \\
      & x_1\in \X_1,  \ \ \ &  &\ x_2\in \X_2,\ \ \ & &\ x_3\in \X_3.
\end{aligned}
\end{equation}
We can write problem~\eqnok{3-stage prob} in a more compact form by using value functions
as discussed in \eqnok{dynpr-1}-\eqnok{dynpr-3}.
More specifically, let $Q_3(x_{2},\xi_{[3]})$ be the 
value function at the third stage and $\Q_3(x_2)$ be the corresponding expected value function
conditionally on $\xi_{[2]}$:
\begin{equation}\label{Defi_sto_V1}
 \begin{array}{lll}
Q_3(x_{2},\xi_{[3]}) &:=&  \min \ f_3(x_3,\xi_3)\\
  &&\text{ s.t.} \ \ A_3 x_3 + B_3 x_2 = b_3, \\
  & &\quad \quad \quad x_3\in \X_3.\\
\Q_3(x_2,\xi_{[2]}) &:=& \bbe_{|\xi_{[2]}} [Q_3(x_{2},\xi_{[3]})].
\end{array}
\end{equation}
We can then define the stochastic cost-to-go function $Q_2(x_{1},\xi_{[2]})$ and its
corresponding (expected) value function as
\begin{equation}\label{Defi_sto_V}
\begin{array}{lll}
Q_2(x_{1},\xi_{[2]}) &:= & \min \  f_2(x_2,\xi_2)+ \Q_{3}(x_2, \xi_{[2]}) \\
 &&\ \text{  s.t.} \ \ A_2 x_2 + B_2x_{1} = b_2,\\
 & & \quad \quad \quad x_2\in \X_2.\\
\Q_2(x_1,\xi_1) &:=& \bbe_{|\xi_1} [Q_2(x_{1},\xi_{[2]}) ] = \bbe [Q_2(x_{1},\xi_2)] .
 \end{array}
 \end{equation}
Problem \eqref{3-stage prob} can then be formulated equivalently as
\begin{equation} \label{first-stage-cp}
\begin{array}{ll}
\min \  f_1(x_1,\xi_1)+ \Q_{2}(x_1,\xi_1) \\
\text{ s.t.} \ \ A_1 x_1 = b_1,\\
\quad \quad \quad x_1 \in \X_1.
\end{array}
\end{equation}

\subsection{Approximate stochastic subgradients}
In order to solve problem~\eqnok{first-stage-cp} by stochastic approximation type methods, we need to
understand how to compute first-order information of the
value functions $\Q_2$ and $\Q_3$.
Recall that for a given closed convex set $\X \subseteq \bbr^n$ and a closed convex function $\Q: \X \to \bbr$,
$g(x)$ is called an $\epsilon$-subgradient of $\Q$ at $x \in X$ if
\begin{equation} \label{def_eps_sub}
\Q(y) \geq \Q(x) + \langle g(x), y-x\rangle -\epsilon \ \ \forall y \in \X.
\end{equation}
The collection of all such $\epsilon$-subgradients of $\Q$ at $x$ is called
the $\epsilon$-subdeifferential of $\Q$ at $x$,
 denoted by $\partial_\epsilon \Q(x)$.
Since both $\Q_2$ and $\Q_3$
are given in the form of (conditional) expectation, their exact first-order information
is hard to compute. We resort to the computation of a stochastic $\epsilon$-subgradient
of these value functions defined as follows.
Observe that we do not assume stagewise independence throughout this section.

\begin{definition}
$G(x,\xi_{[t]})$ is called a stochastic $\epsilon$-subgradient of the value function $\Q_t(x, \xi_{[t-1]}) = \bbe_{|\xi_{[t-1]}}[Q_t(x, \xi_{[t]})]$
if $G(x,\xi_{[t]})$ is an unbiased estimator of an $\epsilon$-subgradient of $\Q_t(x, \xi_{[t-1]})$ with respect to $x$, i.e.,
\begin{equation}\label{def_g}
\bbe_{| \xi_{[t-1]}} [G(x,\xi_{[t]})] = g(x, \xi_{[t-1]}) \ \ \mbox{and} \ \ g(x,\xi_{[t-1]}) \in \partial_\epsilon \Q_t(x, \xi_{[t-1]}).
\end{equation}
\end{definition}

To compute a stochastic $\epsilon$-subgradient of
$\Q_2$  (resp., $\Q_3$), we have to compute an approximate subgradient
of the corresponding stochastic value function $Q_2(x_1, \xi_{[2]})$ (resp., $Q_3(x_2, \xi_{[3])}$).
To this end, we
further assume that strong Lagrange duality holds
 for the optimization problems defined in \eqnok{Defi_sto_V} (resp.,\eqnok{Defi_sto_V1}) almost surely. In other words, these problems can be formulated
 as saddle point problems:
\begin{align}
Q_2(x_1, \xi_{[2]})&=\max_{y_2}\min_{x_2\in \X_2}\langle b_2 - B_2 x_1-A_2 x_2, y_2\rangle + f_2(x_2,\xi_2)+ \Q_3(x_2, \xi_{[2]}),\label{primaldual2}\\
Q_3(x_2, \xi_{[3]})&=\max_{y_3}\min_{x_3\in \X_3}\langle b_3 - B_3 x_2-A_3x_3, y_3\rangle + f_3(x_3,\xi_3).\label{primaldual3}
\end{align}
One set of sufficient conditions to guarantee the equivalence between \eqnok{Defi_sto_V} (resp.,\eqnok{Defi_sto_V1}) and \eqnok{primaldual2} (resp., \eqnok{primaldual3})
is that \eqnok{Defi_sto_V} (resp.,\eqnok{Defi_sto_V1})  is solvable and the Slater condition holds. 

Observe that \eqnok{primaldual2} and \eqnok{primaldual3} are special cases of
a more generic saddle point problem:
\begin{equation} \label{primaldual2.2}
Q(u,\xi) := \max_{y}\min_{x\in \X}\langle b-Bu-Ax, y\rangle + f(x,\xi)+  \Q(x),
\end{equation}
where $b$, $A$ and $B$ are linear mappings of
the random variable $\xi$.
For example, \eqnok{primaldual3} is a special case of  \eqnok{primaldual2.2} with
$u = x_2$, $y = y_3$, $\xi=\xi_{3}$, $f = f_3$ and $ \Q = 0$.
It is worth noting that the first stage problem can also be viewed as a special case of \eqnok{primaldual2.2},
since \eqnok{first-stage-cp} is equivalent to
\begin{equation} \label{primaldual1}
\max_{y_1} \min_{x_1 \in \X_1} \ \left\{ \langle b_1-A_1x_1, y_1\rangle + f_1(x_1,\xi_1)+ \Q_{2}(x_1,\xi_1) \right\}.
\end{equation}

We now discuss how to relate an approximate solution
of the saddle point problem in \eqnok{primaldual2.2} to an approximate
 subgradient of $Q$.
Let $(x_*,y_*)$ be a pair of optimal solutions of the saddle point problem \eqref{primaldual2.2}, i.e.,
\begin{align}
Q(u,\xi) &= \langle y_*,b-Bu-Ax_*\rangle+f(x_*,\xi)+  \Q(x_*)
= f(x_*,\xi)+  \Q(x_*), \label{strong_duality}
\end{align}
where the second identity follows from the complementary slackness of Lagrange duality.
It  can be shown   (see Lemma~1 of \cite{LanZhou17-1}) that
if
\begin{equation}\label{goal_Q}
\begin{aligned}
\gap(\bar z; x,y_*) &:=  \langle y_*, b-Bu-A\bar x\rangle + f(\bar x,\xi)+  \Q(\bar x) \\
& \quad -\langle \bar y, b-Bu-Ax\rangle - f(x,\xi)-  \Q(x) \leq \epsilon,\ \forall x\in \X,
\end{aligned}
\end{equation}
for a given $\bar z:= (\bar x, \bar y) \in Z$ and $u \in \bbr^{n_0}$,
then $-B^T\bar y$ is an $\epsilon$-subgradient of $Q(u, \xi)$ at $u$.

In view of this observation,
in order to compute a stochastic subgradient of $\Q_t(u, \xi_{[t-1]})=\bbe[Q_t(u,\xi_{[t]})|\xi_{[t-1]}]$ at a given point $u$,
we can first generate a random realization $\xi_t$ conditionally on $\xi_{[t-1]}$ and then try to find a pair of solutions $(\bar x, \bar y)$
satisfying 
\begin{align*}
&\langle y_{t*}, b_t-B_tu-A_t \bar x\rangle + f(\bar x,\xi_t)+ \Q_{t+1}(\bar x, \xi_{[t]}) \\
 &-\langle \bar y, b_t -B_tu-A_tx\rangle - f(x,\xi_t)- \Q_{t+1}(x, \xi_{[t]}) \leq \epsilon,\ \forall x\in \X_t,
\end{align*}
where $y_{t*} \equiv y_{t*}(\xi_{[t]})$ denotes the optimal solution for the $t$-th stage problem
associated with
the random realization $\xi_{[t]}$.
We will then use $-B^T \bar y$ as a stochastic $\epsilon$-subgradient
of $\Q_t(u, \xi_{[t-1]})$ at $u$. The difficulty 
exists in that the function $\Q_{t+1}(\bar x, \xi_{[t]})$ is also given in the form
of expectation, and requires a numerical procedure to estimate its value.
We will discuss this in more details in the next section.

\subsection{The DSA algorithm and its convergence properties}

Our goal in this section is to present the basic scheme of the dynamic stochastic
approximation algorithm applied to problem~\eqnok{first-stage-cp}.
This algorithm relies on the following three key primal-dual steps, referred to as stochastic primal-dual
transformation (SPDT), applied to the generic saddle point problem in \eqnok{primaldual2.2}
at every stage.

$(p_+, d_+, \tilde d) = {\rm SPDT}(p, d, d_{\_},  \Q', u, \xi, f, X, \theta, \tau, \eta)$:
\begin{align}
\tilde d &= \theta (d - d_{\_})+ d. \label{def_dual_extra}\\
p_+ &= \argmin_{x \in X} \langle b-Bu-Ax, \tilde d  \rangle + f(x,\xi) +\langle \Q',x \rangle + \tfrac{\tau}{2}\|x - p\|^2. \label{def_primal_proj}\\
d_+ &=  \argmin_{y} \langle -b+Bu+A p_+, y\rangle + \tfrac{\eta}{2} \|y - d\|^2. \label{def_dual_proj}
\end{align}
In the above primal-dual transformation, the input $(p,d, d_{\_})$ denotes the current primal solution, dual solution, and
the previous dual solution, respectively. Moreover, the input $\Q'$ denotes a stochastic
$\epsilon$-subgradient for $\Q$ at the current search point $p$. The parameters $(u, \xi, f, X)$
describe the problem in \eqnok{primaldual2.2} and $(\theta, \tau, \eta)$ are certain algorithmic parameters
to be specified. Given these input parameters, the relation in \eqnok{def_dual_extra} defines a dual extrapolation (or prediction) step
to estimate the dual variable $\tilde d$ for the next iterate. Based on this estimate, \eqnok{def_primal_proj}
performs a primal projection to compute $p_+$, and then \eqnok{def_dual_proj} updates in the dual
space to compute $d_+$ by using the updated $p_+$.
Note that the Euclidean projection in \eqnok{def_primal_proj}
can be extended to the non-Euclidean setting by replacing the term $\|x - p\|^2/2$ with Bregman divergence.
We assume that the above SPDT operator
can be performed very fast or even has explicit expressions.
The primal-dual transformation is closely related to the
alternating direction method of multipliers and the primal-dual hybrid gradient method (see \cite{LanZhou17-1}
for an account of history for this procedure).


In order to solve problem~\eqnok{first-stage-cp}, we will combine the above primal-dual transformation
applied to all the three stages,
the scenario generation for the random variables $\xi_2$ and $\xi_3$
in the second and third stage, and certain averaging steps in both the primal and dual spaces
to compute an approximate pair of primal and dual solution for the saddle point problem in the form of
\eqnok{primaldual2.2} at each stage. More specifically,
the DSA method consists of three loops. The innermost (third) loop runs $N_3$ steps of SPDT in order to compute an approximate stochastic
subgradient 
of the value function $\Q_3$ of the third stage.
The second loop consists of $N_2$ SPDTs applied to the saddle point formulation of the second-stage problem,
which requires the output from the third loop. The outer loop applies $N_1$ SPDTs to the saddle point formulation
of the first-stage optimization problem in \eqnok{first-stage-cp}, using the approximate stochastic  subgradients 
for $\Q_2$ computed by the second loop. In this algorithm, we need to generate $N_1$ and $N_1 \times N_2$ realizations
for the random vectors $\xi_2$ and $\xi_3$, respectively.

Observe that the DSA algorithm
described above is conceptual only since we have not specified any algorithmic parameters $(\theta, \tau, \eta)$ yet.
Two sets of parameters will be chosen depending on the stages where SPDTs are applied.
One set of parameters that may lead to slightly slower rate of convergence but can guarantee the boundedness of
the generated dual iterates will be used for the second and third stage, while another set of
parameter setting is used in the first stage to achieve faster rate of convergence.
Suppose that $f_t$, $t =1, 2, 3$, are generally (not necessarily strongly) convex.
By properly specifying the algorithmic parameters, it can be shown that by setting $N_3 = {\cal O}(1/\epsilon^2)$ and
$N_2 = {\cal O}(1/\epsilon^2)$, we will find an approximate $\epsilon$-solution,
i.e., a point $\bar x_1 \in \X_1$  such that
\[
  \begin{array}{l}
  \bbe[f_1(\bar x_1, \xi_1) + \Q_2(\bar x_1,\xi_1) - (f_1(x_*, \xi_1) + \Q_2(x_*,\xi_1))]  \le \epsilon,\\
\bbe [\| A \bar x_1 - b\|] \le \epsilon,
\end{array}
\]
in at most $N_1 = {\cal O}(1/\epsilon^2)$ outer iterations. As a consequence, the number of
random samples $\xi^2$ and $\xi^3$ used in the DSA method are bounded by
\begin{equation} \label{def_samp_complex}
N_1 = \mathcal{O}(1/\epsilon^2) \ \ \mbox{and} \ \ \ N_1\times N_2 = \mathcal{O}(1/\epsilon^4),
\end{equation}
respectively.
Now consider the case when  the objective functions $f_t$, $t = 1, 2, 3$, are strongly convex.
In that case the above complexity results can be significantly improved. In particular,
by choosing $N_3 = {\cal O}(1/\sqrt{\epsilon})$ and
$N_2 = {\cal O}(1/\epsilon)$, we will find an approximate $\epsilon$-solution
in at most $N_1 = {\cal O}(1/\epsilon)$ outer iterations. As a result, it  can be shown  that
the number of random samples of  $\xi_2$ and $\xi_3$ will be bounded by $N_1$ and $N_1\times N_2$, i.e., $\mathcal{O}(1/\epsilon)$ and $\mathcal{O}(1/\epsilon^2)$, respectively,

It should be noted that our analysis of DSA focuses on the optimality of the first-stage decisions, and
the decisions we generated for the later stages are mainly used for computing the approximate stochastic
subgradients for the value functions at each stage.
Except for the first stage decision $\bar x_1$,
the performance guarantees (e.g., feasibility and optimality)
 that we can provide for later stages are
dependent on the sequences of random variables (or scenarios) we generated.

\subsection{DSA for general multistage stochastic optimization}

In this section, we consider a multistage stochastic optimization problem given by
\begin{equation} \label{multi-stage prob}
\begin{array}{ll}
\min \  f_1(x_1,\xi_1)+ \Q_{2}(x_1,\xi_1)  \\
\text{ s.t.} \ \ A_1 x_1 = b_1,\\
\quad \quad \quad x_1 \in \X_1,
\end{array}
\end{equation}
where the value factions $\Q_t$, $t =2, \ldots, T$, are recursively defined by
\begin{equation}\label{Defi_sto_V_m}
\begin{array}{lll}
 \Q_t(x_{t-1},\xi_{[t-1]}) &:=& F_{t-1} (x_{t-1},\xi_{t-1}) + \bbe_{|\xi_{[t-1]}} [Q_t(x_{t-1},\xi_{[t]})], \ \ t = 2, \ldots, T-1,\\
Q_t(x_{t-1},\xi_{[t]}) &:= & \min \ f_t(x_t,\xi_t)+ \Q_{t+1}(x_t)\\
 &&\ \text{  s.t.} \ \ A_t x_t + B_t x_{t-1} = b_t,\\
 & & \quad \quad \quad x_t \in \X_t,
 \end{array}
 \end{equation}
 and
\begin{equation}\label{Defi_sto_V1_m}
 \begin{array}{lll}
  \Q_T(x_{T-1},\xi_{[T-1]}) &:=& \bbe_{\xi_T|\xi_{[T-1]}} [Q_T(x_{T-1},\xi_{[T]})],\\
Q^T(x_{T-1},\xi_{[T]}) &:=&  \min \ f_T(x_T,\xi_T)\\
  &&\text{ s.t.} \ \ A_T x_T  + B_T x_{T-1} = b_T, \\
  & &\quad \quad \quad x_T\in \X_T.
\end{array}
\end{equation}
Here $\xi_t$ are random variables, $f_t(\cdot,\xi_t)$ are relatively simple functions, and $F_t(\cdot,\xi_t)$ are general (not necessarily simple)
Lipschitz continuous convex functions. We also assume that one can compute the subgradient $F'(x_t,\xi_t)$ of
function $F_t(x_t,\xi_t)$ at any point $x_t\in \X_t$ for a given $\xi_t$.

Problem~\eqnok{multi-stage prob} is more general than problem \eqref{3-stage prob} (or equivalently problem~\eqnok{first-stage-cp}) in the
following sense. First, we are dealing with a more complicated multistage stochastic optimization problem where the number of stages $T$ in \eqnok{multi-stage prob} can be greater than three.
Second, the value function $ \Q_t(x_{t-1}, \xi_{[t-1]})$ in \eqnok{Defi_sto_V_m}
is defined as the summation of  $F_{t-1} (x_{t-1},\xi_{t-1})$ and $ \bbe_{|\xi_{[t-1]}}[Q_t(x_{t-1},\xi_{[t]})]$,
where $F_{t-1}$ is not necessarily simple.
The DSA algorithm can be generalized for solving problem~\eqnok{multi-stage prob}.
More specifically, we will call the SPDT operators to compute a stochastic $\epsilon$-subgradient of $\Q_{t+1}$
at $x_t$, $t= 1, \ldots, T-2$, in a recursive manner until we obtain the $\epsilon$-subgradient of $\Q_T$ at $x_{T-1}$.

Similar to the three-stage problem, let $N_t$ be the number of iterations for stage $t$ subproblem.
For the last stage $T$, we will set $N_T = {\cal O}(1/\epsilon)$ and
$N_T = {\cal O} (1/\sqrt{\epsilon})$ for the generally convex and strongly convex cases, respectively.
For the middle stages $t = 2, \ldots, T-1$, we set $N_t = {\cal O} (1/\epsilon^2)$ and $N_t = {\cal O}(1/\epsilon)$.
Different algorithmic settings will be used for either generally convex or strongly convex cases.
Moreover, less aggressive stepsize selection will be used for the inner loops to guarantee the boundedness of
the generated dual variables. Under these parameter selection, we can show that
the number of outer loops performed by the DSA method to find an approximate $\epsilon$-solution
can be bounded by $N_1 = {\cal O}(1/\epsilon^2)$ and ${\cal O}(1/\epsilon)$, respectively, for the
generally convex and strongly convex problems.

In view of these results,
the total number of scenarios required to find an $\epsilon$-solution of \eqnok{multi-stage prob} is given by $N_2 \times N_3 \times \ldots N_T$,
and hence will grow exponentially with respect to $T$, no matter whether the objective functions are strongly convex or not.
These sampling complexity bounds match well with
those lower bounds in \cite{ShaNem04,sha06}, implying that
multi-stage stochastic optimization problems are essentially intractable for $T \ge 5$ and a moderate target accuracy.
Hence, it is reasonable to use the DSA algorithm only for multi-stage stochastic optimization problems with $T=3$ or $4$ and $\epsilon$ relatively large.
However, it is interesting to point out that the DSA algorithm
only needs to go through the scenario tree once and hence its memory requirement increases only linearly with respect to $T$.


\subsection{Combined EDDP and DSA for hierarchical problems}
In this section, we discuss how EDDP/SDDP type algorithms can be applied
together with DSA type algorithms for solving a class of
hierarchical stationary stochastic programs (HSSPs).
 HSSPs can model problems with a hierarchy of decision-making, e.g., how managerial decisions influence day-to-day operations in a factory.
 More specifically, the upper level decisions are made by solving a stationary stochastic program over an infinite horizon (see Section~\ref{sec-inf}),
where the number of stages $T = \infty$, the cost function $f_t = f$,
the feasible set $\X_t = \X$, and the random variables $\xi_t$ are iid having
the same distribution as $\xi$.
The lower level decisions  are determined by a stochastic two-stage program,
\begin{align} \label{eq:qgp_a11}
  f(x,\xi) := \min_{z_1 \in Z_1(x,\xi)} f_1(z_1,\xi) + \mathbb{E}_{\zeta}
  \big[ \min_{z_2 \in Z_2(z_1,\zeta)} f_2(z_2,\zeta) \big],
\end{align}
where $\zeta$ denotes the random variable  independent of $\xi$ involved in the second stage problem.

We extend the dual dynamic programming method to a so-called \textit{hierarchical dual dynamic programming} by accounting for inexact solutions to the subproblems. To solve the lower-level stochastic multi-stage problem, we approximate it using SAA approach and solve the resulting problem using a primal-dual stochastic approximation method. This method can be generalized to a DSA-type method when the number of stages is three or more. Our results show that when solving an infinite-horizon hierarchical problem where the top-level decision variable is of modest size (i.e., $n=4$) and the lower-level stochastic multistage program has a modest number of stages (2 or 3), then by integrating dual dynamic programming to handle the top-level decisions and stochastic approximation-type method for the lower-level decisions, we can get a polynomial computational complexity that is independent of the dimension of the lower-level decision variables (see \cite{JuLan23} for details).


\bibliographystyle{plain}
\bibliography{references}

\begin{thebibliography}{10}

\bibitem{ADEH:1999}
P.~Artzner, F.~Delbaen, J.-M. Eber, and D.~Heath.
\newblock Coherent measures of risk.
\newblock {\em Mathematical Finance}, 9:203--228, 1999.

\bibitem{beal:55}
E.~M.~L. Beale.
\newblock On minimizing a convex function subject to linear inequalities.
\newblock {\em Journal of the Royal Statistical Society, Series B},
  17:173--184, 1955.

\bibitem{Bel}
R.E. Bellman.
\newblock {\em Dynamic Programming}.
\newblock Princeton University Press, 1957.

\bibitem{ber78}
D.P. Bertsekas and S.E. Shreve.
\newblock {\em Stochastic Optimal Control, The Discrete Time Case}.
\newblock Academic Press, New York, 1978.

\bibitem{bir85}
J.R. Birge.
\newblock Decomposition and partitioning methods for multistage stochastic
  linear programs.
\newblock {\em Operations Research}, 33:989--1007, 1985.

\bibitem{birge2011}
R.~Birge and F.~Louveaux.
\newblock {\em Introduction to Stochastic Programming}.
\newblock Springer, New York, 2nd edition, 2011.

\bibitem{BS2000}
J.~F. Bonnans and A.~Shapiro.
\newblock {\em Perturbation Analysis of Optimization Problems}.
\newblock Springer Series in Operations Research. Springer, 2000.

\bibitem{dant:55}
G.B. Dantzig.
\newblock Linear programming under uncertainty.
\newblock {\em Management Science}, 1:197--206, 1955.

\bibitem{ding2019}
L.~Ding, S.~Ahmed, and A.~Shapiro.
\newblock A python package for multi-stage stochastic programming.
\newblock {\em Optimization online}, 2019.

\bibitem{dyer}
M.~Dyer and L.~Stougie.
\newblock Computational complexity of stochastic programming problems.
\newblock {\em Mathematical Programming}, 106:423--432, 2006.

\bibitem{follm}
H.~F\"ollmer and A.~Schied.
\newblock {\em Stochastic Finance: An Introduction in Discrete Time}.
\newblock Walter de Gruyter, Berlin, 2nd edition, 2004.

\bibitem{FR2021}
C.~F\"ullner and S.~Rebennack.
\newblock Stochastic dual dynamic programming and its variants.
\newblock {\em Optimization online}, 2021.

\bibitem{Gir:2014}
P.~Girardeau, V.~Lecl\'ere, and A.~B. Philpott.
\newblock On the convergence of decomposition methods for multistage stochastic
  convex programs.
\newblock {\em Mathematics of Operations Research}, 40:130--145, 2016.

\bibitem{GSC2022}
V.~Guigues, A.~Shapiro, and Y.~Cheng.
\newblock Risk-averse stochastic optimal control: an efficiently computable
  statistical upper bound.
\newblock {\em Optimization online}, 2022.

\bibitem{GSC22}
V.~Guigues, A.~Shapiro, and Y.~Cheng.
\newblock Duality and sensitivity analysis of multistage linear stochastic
  programs.
\newblock {\em European Journal of Operational Research}, 308:752--767, 2023.

\bibitem{HKW}
G.A. Hanasusanto, D.~Kuhn, and W.~Wiesemann.
\newblock A comment on computational complexity of stochastic programming
  problems.
\newblock {\em Mathematical Programming}, 159:557--569, 2015.

\bibitem{JuLan23}
C.~Ju and G.~Lan.
\newblock Dual dynamic programming for stochastic programs over an infinite
  horizon.
\newblock {\em arXiv}, 2023.

\bibitem{kel:60}
J.E. Kelley.
\newblock The cutting-plane method for solving convex programs.
\newblock {\em Journal of the Society for Industrial and Applied Mathematics},
  8:703--712, 1960.

\bibitem{Lan10-3}
G.~Lan.
\newblock An optimal method for stochastic composite optimization.
\newblock {\em Mathematical Programming}, 133(1):365--397, 2012.

\bibitem{LAN2015}
G.~Lan.
\newblock Bundle-level type methods uniformly optimal for smooth and nonsmooth
  convex optimization.
\newblock {\em Mathematical Programming}, 149(1-2):1--45, 2015.

\bibitem{LanBook2020}
G.~Lan.
\newblock {\em First-order and Stochastic Optimization Methods for Machine
  Learning}.
\newblock Springer Nature, Switzerland AG, 2020.

\bibitem{Lan22}
G.~Lan.
\newblock Complexity of stochastic dual dynamic programming.
\newblock {\em Mathematical Programming}, 191:717--754, 2022.

\bibitem{lns11}
G.~Lan, A.~S. Nemirovski, and A.~Shapiro.
\newblock Validation analysis of mirror descent stochastic approximation
  method.
\newblock {\em Mathematical Programming}, 134:425--458, 2012.

\bibitem{LanZhou17-1}
G.~Lan and Z.~Zhou.
\newblock Dynamic stochastic approximation for multi-stage stochastic
  optimization.
\newblock {\em Mathematical Programming}, 187:487--532, 2021.

\bibitem{Leclere}
V.~Lecl\'ere, P.~Carpentier, J-P Chancelier, A.~Lenoir, and F.~Pacaud.
\newblock Exact converging bounds for {S}tochastic {D}ual {D}ynamic
  {P}rogramming via {F}enchel duality.
\newblock {\em SIAM J. Optimization}, 30:1223--1250, 2020.

\bibitem{LNN1995}
C.~Lemar\'{e}chal, A.~Nemirovskii, and Y.~Nesterov.
\newblock New variants of bundle methods.
\newblock {\em Mathematical Programming}, 69:111--147, 1995.

\bibitem{Lohndorf}
N.~L\"ohndorf and A.~Shapiro.
\newblock Modeling time-dependent randomness in stochastic dual dynamic
  programming.
\newblock {\em European Journal of Operational Research}, 273:650--661, 2019.

\bibitem{NJLS09-1}
A.~S. Nemirovski, A.~Juditsky, G.~Lan, and A.~Shapiro.
\newblock Robust stochastic approximation approach to stochastic programming.
\newblock 19:1574--1609, 2009.

\bibitem{nemyud:83}
A.~S. Nemirovski and D.~Yudin.
\newblock {\em Problem complexity and method efficiency in optimization}.
\newblock Wiley-Interscience Series in Discrete Mathematics. John Wiley, XV,
  1983.

\bibitem{Nesnem}
Y.~Nesterov and A.~Nemirovski.
\newblock {\em Interior-Point Polynomial Algorithms in Convex Programming}.
\newblock SIAM, Philadelphia, 1994.

\bibitem{per1991}
M.V.F. Pereira and L.M.V.G. Pinto.
\newblock Multi-stage stochastic optimization applied to energy planning.
\newblock {\em Mathematical programming}, 52(1-3):359--375, 1991.

\bibitem{pflug2000}
G.~Pflug.
\newblock Some remarks on the value-at-risk and the conditional value-at-risk.
\newblock In {\em Probabilistic Constrained Optimization: Methodology and
  Applications, S. Uryasev (Ed.)}. Kluwer Academic Publishers, Norwell, MA,
  2000.

\bibitem{phil13}
A.B. Philpott, V.L. de~Matos, and E.~Finardi.
\newblock On solving multistage stochastic programs with coherent risk
  measures.
\newblock {\em Operations Research}, 61(4):957--970, 2013.

\bibitem{pichsha}
A.~Pichler and A.~Shapiro.
\newblock Mathematical foundations of distributionally robust multistage
  optimization.
\newblock {\em SIAM J. Optimization}, 31:3044--3067, 2021.

\bibitem{pol90}
B.T. Polyak.
\newblock New stochastic approximation type procedures.
\newblock {\em Automat. i Telemekh.}, 7:98--107, 1990.

\bibitem{pol92}
B.T. Polyak and A.B. Juditsky.
\newblock Acceleration of stochastic approximation by averaging.
\newblock {\em SIAM J. Control and Optimization}, 30:838--855, 1992.

\bibitem{powell}
W.B. Powell.
\newblock {\em Approximate Dynamic Programming: Solving the Curses of
  Dimensionality}.
\newblock John Wiley and Sons, New York, 2nd edition, 2011.

\bibitem{RobMon51-1}
H.~Robbins and S.~Monro.
\newblock A stochastic approximation method.
\newblock {\em Annals of Mathematical Statistics}, 22:400--407, 1951.

\bibitem{roc74}
R.~T Rockafellar.
\newblock {\em Conjugate Duality and Optimization}.
\newblock Society for Industrial and Applied Mathematics, Philadelphia, 1974.

\bibitem{roc1970}
R.T. Rockafellar.
\newblock {\em Convex Analysis}.
\newblock Princeton University Press, 1970.

\bibitem{ury2}
R.T Rockafellar and S.~Uryasev.
\newblock Conditional value-at-risk for general loss distributions.
\newblock {\em J. of Banking and Finance}, 26(7):1443--1471, 2002.

\bibitem{RuSh:2006b}
A.~Ruszczy\'{n}ski and A.~Shapiro.
\newblock Conditional risk mappings.
\newblock {\em Mathematics of Operations Research}, 31:544--561, 2006.

\bibitem{RuSh:2006a}
A.~Ruszczy\'{n}ski and A.~Shapiro.
\newblock Optimization of convex risk functions.
\newblock {\em Mathematics of Operations Research}, 31:433--452, 2006.

\bibitem{shaejor2011}
A.~Shapiro.
\newblock Analysis of stochastic dual dynamic programming method.
\newblock {\em European Journal of Operational Research}, 209:63--72, 2011.

\bibitem{shacheng2021}
A.~Shapiro and Y.~Cheng.
\newblock Central limit theorem and sample complexity of stationary stochastic
  programs.
\newblock {\em Operations Research Letters}, 49:676--681, 2021.

\bibitem{shacheng2021a}
A.~Shapiro and Y.~Cheng.
\newblock Dual bounds for periodical stochastic programs.
\newblock {\em Operations Research}, 2021.

\bibitem{SDR}
A.~Shapiro, D.~Dentcheva, and A.~Ruszczy\'{n}ski.
\newblock {\em Lectures on Stochastic Programming: Modeling and Theory}.
\newblock SIAM, Philadelphia, third edition, 2021.

\bibitem{shading}
A.~Shapiro and L.~Ding.
\newblock Periodical multistage stochastic programs.
\newblock {\em SIAM J. Optimization}, 30:2083--2102, 2020.

\bibitem{ShaNem04}
A.~Shapiro and A.~Nemirovski.
\newblock On complexity of stochastic programming problems.
\newblock {\em E-print available at: http://www.optimization-online.org}, 2004.

\bibitem{shanem}
A.~Shapiro and A.~Nemirovski.
\newblock On complexity of stochastic programming problems.
\newblock In V.~Jeyakumar and A.M. Rubinov, editors, {\em Continuous
  Optimization: Current Trends and Applications: Current Trends and
  Applications}, pages 111--144. Springer, 2005.

\bibitem{Sha2012a}
A.~Shapiro, W.~Tekaya, J.P. da~Costa, and M.~Pereira Soares.
\newblock Risk neutral and risk averse stochastic dual dynamic programming
  method.
\newblock {\em European Journal of Operational Research}, 224:375--391, 2013.

\bibitem{sha06}
A.~Shaprio.
\newblock On complexity of multistage stochastic programs.
\newblock {\em Operations Research Letters}, 34:1--8, 2006.

\bibitem{ZhangSun22-1}
S.~Zhang and X.~Sun.
\newblock Stochastic dual dynamic programming for multistage stochastic
  mixed-integer nonlinear optimization.
\newblock {\em Mathematical Programming}, 196:935–985, 2022.

\bibitem{zipkin}
P.H. Zipkin.
\newblock {\em Foundation of inventory management}.
\newblock McGraw-Hill, 2000.

\end{thebibliography}

\end{document}